\documentclass{article}
\usepackage[fontsize=13pt]{scrextend}
\usepackage{amsmath}
\usepackage{amsthm}
\usepackage{amssymb}
\usepackage{graphicx} 
\usepackage{multicol} 
\usepackage{enumerate} 
\usepackage[british]{babel}
\usepackage{tikz-cd} 
\usepackage{mathrsfs} 
\usepackage[hidelinks]{hyperref}
\usepackage{multicol}

\usepackage{lipsum}

\theoremstyle{plain}  
\newtheorem{theorem}{Theorem}[section]
\newtheorem{proposition}[theorem]{Proposition}
\newtheorem{corollary}[theorem]{Corollary}
\newtheorem{lemma}[theorem]{Lemma}
\newtheorem{conjecture}[theorem]{Conjecture}

\theoremstyle{definition} 
\newtheorem{definition}[theorem]{Definition}
\newtheorem{example}[theorem]{Example}

\theoremstyle{remark} 
\newtheorem{remark}[theorem]{Remark}

\newcommand{\Q}{\mathbb{Q}}  
\newcommand{\N}{\mathbb{N}}  
\newcommand{\C}{\mathbb{C}}  
\newcommand{\A}{\mathscr{A}} 
\newcommand{\K}{\mathscr{K}} 
\newcommand{\icis}{\textsc{icis} } 
\newcommand{\singular}{\textsc{Singular} } 
\newcommand{\codimAe}{\textup{codim}_{\A_e}} 
\newcommand{\codimKe}{\textup{codim}_{\K_e}} 
\newcommand{\depth}{\textup{depth}\,} 
\newcommand{\grade}{\textup{grade}\,} 
\newcommand{\Ann}{\textup{Ann}\,} 
\newcommand{\pd}{\textup{pd}} 
\newcommand{\height}{\textup{ht}\,} 
\renewcommand{\O}{\mathscr{O}} 
\newcommand{\F}{\mathscr{F}} 
\newcommand{\X}{\mathfrak{X}} 
\newcommand{\M}{\mathscr{M}} 
\newcommand{\m}{\mathfrak{m}} 
\newcommand{\q}{\mathfrak{q}} 
\newcommand{\length}{\textup{length}\,} 
\newcommand{\id}{\textup{id}} 
\newcommand{\im}{\textup{im}\,} 
\newcommand{\Tor}{\textup{Tor}\,} 
\newcommand{\rel}{\textup{rel}} 
\renewcommand{\sp}{\textup{sp}} 

\renewcommand{\phi}{\varphi} 
\newcommand{\gr}{\text{gr}} 
\newcommand{\Derlog}{\textup{Derlog}\,} 
\newcommand{\Hom}{\textup{Hom}} 
\newcommand{\Ext}{\textup{Ext}} 

\usepackage{txfonts} \usepackage[T1]{fontenc}

\usepackage[ backend=biber, style=alphabetic, sorting=nyvt]{biblatex}
\bibstyle{alpha} 
\usepackage{csquotes}

\usepackage{scalefnt}
\let\oldtextsc\textsc
\renewcommand{\textsc}[1]{\oldtextsc{\scalefont{1.1}#1}}

\author{Alberto Fernández Hernández}
\date{}

\begin{document}
\pagenumbering{roman}
\begin{titlepage}
    \vspace{0.5in}
    \begin{center}
    \begin{figure}[ht]
    \begin{center}
    \includegraphics[width=0.49\textwidth]{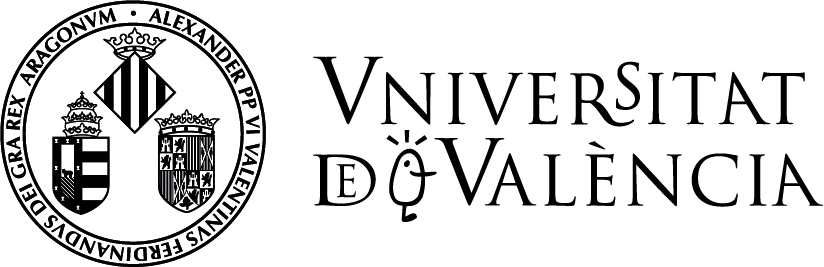}
    \includegraphics[width=0.49\textwidth]{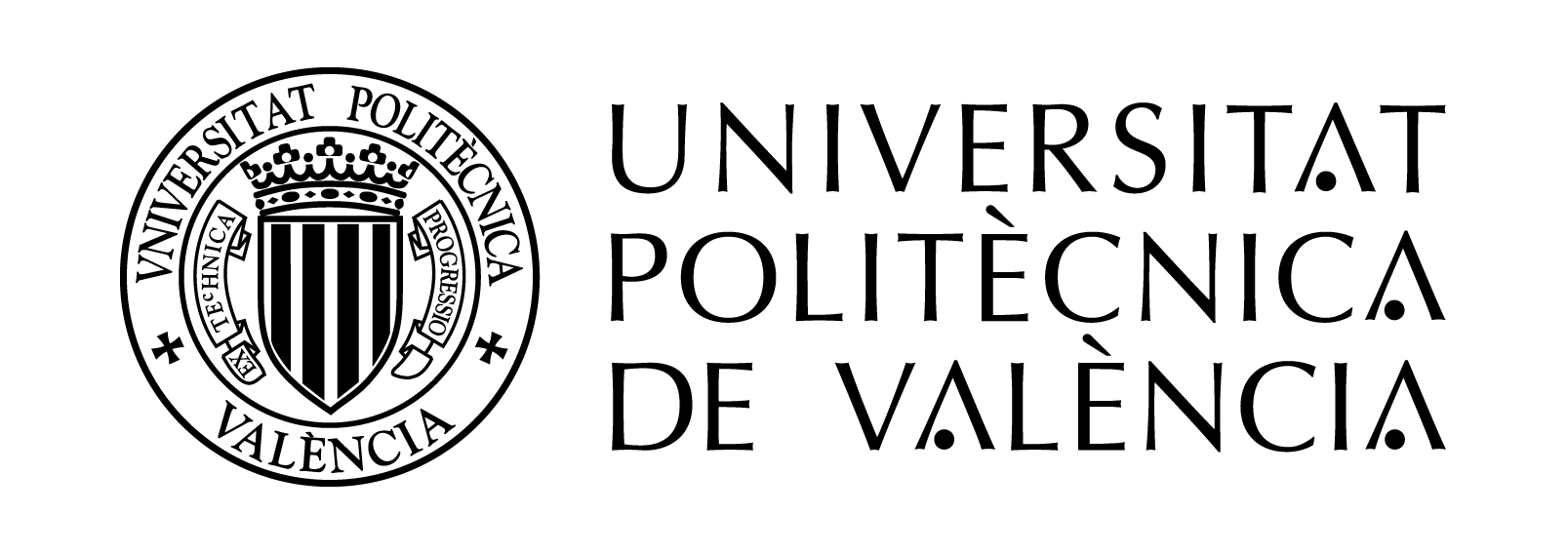}
    \end{center}
    \end{figure}
    \vspace{0.5in}
    \begin{huge}
    Master's Final Thesis - 2022/2023

\vfill    
    \textbf{Disentangling mappings \\ defined on \icis} \\
    \end{huge}
    
\vfill    
    \begin{Large}
    Author: {\bf
    Alberto Fern\'andez Hern\'andez} \\
    \end{Large}
    
    \vspace*{0.2in}
    
    \begin{Large}
    Advisor: {\bf \sc
    Juan Jos\'e Nu\~{n}o Ballesteros } \\
    \end{Large}
    \vspace{0.5in}
    \rule{110mm}{0.3mm}
    \vspace{0.3in}
    \begin{figure}[ht]
        \centering
        \includegraphics[width=0.3\textwidth]{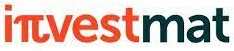}
    \end{figure}

    \begin{Large}
        Máster en Investigación Matemática
    \end{Large}
    \end{center}
    \end{titlepage}
    
\newpage
\textcolor{white}{blank page}
\thispagestyle{empty}

\newpage
\thispagestyle{empty}

\begin{flushright}
    ``To know that we know what we know, and to know that we \\ do not know what we do not know, that is true knowledge'' \\ Nicolaus Copernicus
\end{flushright}
\vfill 
\begin{center} \begin{Large}\textbf{Abstract} \end{Large}\end{center}

Let $(X,S)$ be an isolated complete intersection singularity of dimension $n$, and let $f:(X,S)\rightarrow (\C^{n+1},0)$ be a germ of $\A$-finite mapping. In this project, our main contribution is that we show the case $n=2$ of the general Mond conjecture, which states that $\mu_I(X,f)\geq \codimAe(X,f)$, with equality provided $(X,f)$ is weighted homogeneous. Before this project, the only known case for which the conjecture was known to hold is in the case that $n=1$ and $(X,S)$ is a plane curve. In order to prove the case $n=2$ of the generalised Mond conjecture, we extend the constructions of \cite{BobadillaNunoPenafort} to this general framework. 

\vfill
\begin{figure}[h]
    \centering
    \includegraphics[scale=0.6]{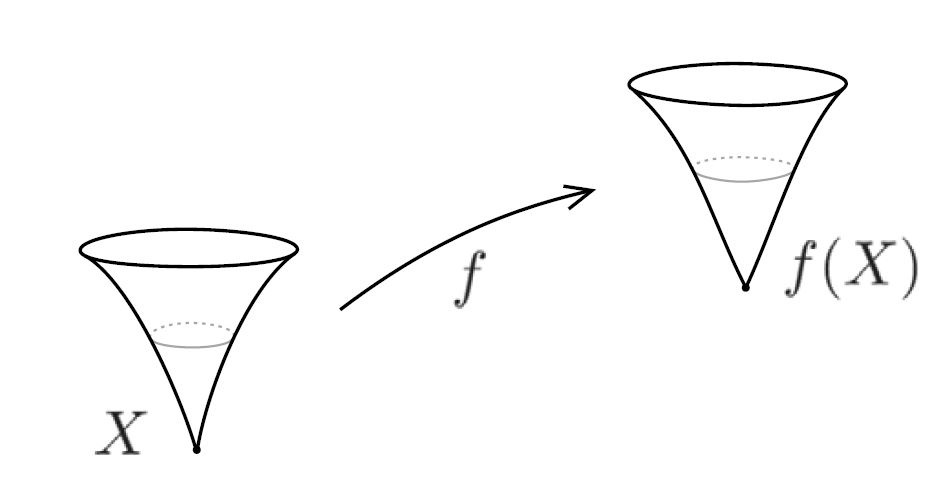}
    \caption{Representation of a mapping $f$ defined on an \icis $X$, and its image $f(X)$. Image extracted from the
article \cite{roberto} by R. Giménez Conejero and J.J. Nuño-Ballesteros.} 
\end{figure}
\vfill
\newpage
\thispagestyle{empty}
\textcolor{white}{blank page}
\newpage
\thispagestyle{empty}
\tableofcontents
\newpage

\section*{Introduction}
\markboth{Introduction}{Introduction}
\addcontentsline{toc}{section}{Introduction}
\setcounter{page}{1}
\pagenumbering{arabic}
The main objects that are studied along this project are the hypersurface singularities. As it is known, if a hypersurface has isolated singularity (meaning that the points of the hypersurface where it is singular is isolated), then one can define two different invariants that capture some information of the hypersurface, namely the Milnor and Tjurina numbers, $\mu$ and $\tau$, respectively. The first of them, $\mu$, has a more topological flavour, since it captures the homotopy type of a smooth deformation of the hypersurface (namely, the Milnor fibre). On the other hand, the Tjurina number has a much more rigid nature, and it is related with the number of degrees of freedom that are required to completely deform the hypersurface. An inspection of the algebraic formulas that define them makes obvious that $\mu \geq \tau$, with equality in the case that the hypersurface has a simple form, namely if it is weighted homogeneous. \\

In 1991, David Mond extended the definition of $\mu$ to a particular kind of hypersurfaces which could be parametrised through a mapping $f:(\C^n,S)\rightarrow (\C^{n+1},0)$ with a nice behaviour (meaning that it has isolated instability, or, equivalently, that $f$ is $\A$-finite). He then introduced the concept of \textit{image Milnor number} $\mu_I(f)$ of the parametrisation $f$. This new invariant played the role of the Milnor number for this kind of hypersurfaces that do not have isolated singularity anymore. Moreover, the codimension of $f$, denoted by $\codimAe (f)$, is an invariant that measures the number of degrees of freedom to fully deform $f$, in the same way the Tjurina number does. If $\mu_I(f)$ and $\codimAe(f)$ are the natural versions of $\mu$ and $\tau$ in this different context, the question is therefore evident:
\begin{center}
    Is it true that $\mu_I(f)\geq \codimAe (f)$?
\end{center}
And equality does hold in the weighted homogeneous case, as it happened in the previous setting? This question is currently known as the Mond conjecture, and it is still an open problem in Singularity Theory. The only known cases are $n=1,2$. \\

Our interest is to extend the scope of the Mond conjecture to more hypersurfaces. It turns out that some hypersurfaces do not admit a parametrisation, but they do admit what we call a \textit{normalisation} over a normal space. A possible generalisation is therefore to study hypersurfaces that admit a normalisation over an \icis $(X,S)$, which is a space with isolated singularity and that where the number of equations that defines it is precisely its codimension; and where the normalisation mapping has isolated instability. In this context, the Mond conjecture can be generalised to mappings $f: (X,S)\rightarrow (\C^{n+1},0)$ to yield whether $\mu_I(X,f)\geq \codimAe(X,f)$. This general version of the conjecture is only known to be true whenever $n=1$ and $(X,S)$ is a plane curve. Our main contribution in this project is to show the case $n=2$ of this general conjecture. Therefore, we provide a positive answer for surfaces in the three-dimensional space $\C^3$.\\

In order to do so, we settle first all the required prerequisites along this project. 
In chapter 1, we study the theory of deformations of mappings defined on a smooth source, which is the context in which the classic conjecture is stated. Then, we analyse the geometry and the topology of \textsc{icis} in chapter 2, and we present the definitions of the Milnor and Tjurina numbers. In chapter 3, we state the classic Mond conjecture motivated by the previous chapters. Chapter 4 extends the main definitions and theorems of the first chapter to the framework of mappings defined on \textsc{icis}. In chapter 5, we perform the construction of a key module that provides an algebraic formula to compute the image Milnor number in the classic case. In chapter 6, we provide a generalisation of the results shown in chapter 5 for mappings defined on \textsc{icis}, and we prove the case $n=2$ of the generalised Mond conjecture. Lastly, the appendix contains a motley collection of definitions and results from commutative algebra that are required in the more advanced chapters, and we include them for the sake of completeness. \\

The main reference of chapters 1 and 3 is the complete book \cite{juanjo}, which develops the singularity theory of mappings with smooth source. For chapter 2, we follow a combination of the book \cite{Looijenga} and the clear notes \cite{icis}. Chapter 4 is based on the article \cite{Mond-Montaldi} of David Mond and James Montaldi, where deformations of singularities of mappings defined on \textsc{icis} were first analysed, and we also follow the brief and clear exposition of this topic given in \cite{roberto}. Chapter 5 is a review of the recent article \cite{BobadillaNunoPenafort}, where the module $M(g)$ that helps in measuring the image Milnor number was first introduced. Chapter 6 is our main contribution, and extends the results of the previous article to mappings defined on \textsc{icis} to prove the case $n=2$ of the general Mond conjecture. Lastly, the appendix is a combination of the books \cite{Matsumura}, \cite{Eilenberg} and \cite{juanjo}. 

\newpage
\thispagestyle{empty}
\textcolor{white}{.}

\newpage
\section{Prerequisites on singularity theory of mappings}
Let us assume that $U\subset \C^n$ is an open set containing the origin, and that we have a holomorphic mapping $f:U\subset \C^n\rightarrow \C^p$ which we seek to analyze in terms of its local behavior at the origin. In such a case, the most natural course of action is to compute the jacobian matrix of $f$ at $0$. If this matrix has maximum rank, then the origin is a regular point of $f$, and its local behavior can be determined using the implicit function theorem. If $n\leq p$, then $f$ has an immersive point at the origin, whereas if $n\geq p$, it has a submersive point. Singularity theory comes into play when neither of these cases occur: when the rank of the jacobian matrix is not maximum, the local behavior of the function becomes much more intricate. In this section, we provide a brief overview of the singularity theory of applications developed by Mather, which arose from the need to understand the local behavior of mappings. For a comprehensive reference on this subject, we recommend \cite{juanjo}.\\

The initial step in studying singularities of mappings involves identifying how to differentiate between seemingly different singularities. However, this does not have an absolute answer, and there are several paths that can be taken from this point. Various equivalences between singularities of mappings offer alternative approaches to their study. The subsequent sections introduce the two primary parallel methods of examining singularities of mappings: namely, the $\A$-equivalence and the $\K$-equivalence. Despite the fact that numerous definitions apply to both forms of equivalence, they will be examined individually. It is important to note that these two distinct methods for analyzing singularities of mappings should not be viewed as opposing approaches. In fact, they complement each other quite effectively. $\A$-equivalence is more relevant, but also more sophisticated, whereas $\K$-equivalence is somewhat simpler in certain regards, yet provides valuable insights that supplement the former method in many respects.\\

We will assume that the reader possesses some familiarity with germ notation, where $(\C^n,S)$ denotes the multigerm of $\C^n$ in proximity to the points of $S=\{s_1, \ldots, s_m\}\subset \C^n$. Essentially, this means that our focus lies on the behavior of mappings close to the points of $S$. In the vast majority of examples, we will set $S={0}$, allowing us to analyze the behavior of mappings near the origin. However, we will develop the theory in this general context, as multigerms are valuable in certain specific cases.
\subsection{$\A$-equivalence of mappings}
Within this section, we examine the $\A$-equivalence of map-germs, which materialize as a result of changes of coordinates in both the source and target spaces. The definition of $\A$-equivalence of map-germs is, in a broad sense, the most intuitive manner in which to differentiate between different singularities of mappings. It seems reasonable that, if a map-germ is altered by modifications to the source and target coordinates, the behavior of the perturbed mapping would be equivalent to that of the original one.
\begin{definition} Let $f,g:(\C^n,S)\rightarrow (\C^p,0)$ be germs of holomorphic mappings. We say that $f$ and $g$ are $\A$\textit{-equivalent} if there exist isomorphisms (\textit{i.e.}, biholomorphisms) $\phi:(\C^n,S)\rightarrow (\C^n,S)$ and $\psi:(\C^p,0)\rightarrow (\C^p,0)$ such that $g=\psi\circ f \circ \phi^{-1}$. In other words, we ask the diagram
    \begin{center}
        \begin{tikzcd}
          (\C^n,S) \arrow[r, "f"] \arrow[d, "\phi"]
          & (\C^p,0) \arrow[d, "\psi"] \\
          (\C^n,S) \arrow[r, "g"]
          & (\C^p,0)
          \end{tikzcd}
      \end{center}
    to commute. 
\end{definition}
The notion of $\A$-equivalence defines an equivalence relation in the set of map-germs, which is straightforward to verify. The ultimate goal of singularity theory of mappings is to provide a comprehensive classification of the map-germs $f:(\C^n,S)\rightarrow (\C^p,0)$ under the $\A$-equivalence relation. However, this task becomes increasingly challenging as the dimensions $(n,p)$ grow.
\begin{remark}\label{remark:Aequivalence} It is usual to consider the $\A$-equivalence through an action in the space of mappings. More formally, if we define $\A=\text{Isom}(\C^n,S)\times \text{Isom}(\C^p,0)$ as the set of pairs of isomorphisms $(\phi, \psi)$, one can define a natural action of $\A$ on the set of germs $f:(\C^n,S)\rightarrow (\C^p,0)$ by 
\begin{equation*}
    (\phi, \psi)\cdot f = \psi \circ f \circ \phi^{-1}.
\end{equation*}
Hence, the germs $f,g:(\C^n,S)\rightarrow (\C^p,0)$ are $\A$-equivalent if and only if there exists a pair $(\phi, \psi)\in \A$ such that $g=(\phi, \psi)\cdot f$. In other words, the mappings that are $\A$-equivalent to $f$ are precisely the elements of the orbit of $f$ though the action induced by $\A$. 
\end{remark}
One of the main concepts in the singularity theory of mappings is stability. Roughly speaking, one says that a map-germ $f$ is stable if every \textit{deformation} $f_t$ of $f$ has the same kind of singularity as $f$, that is, one can write $f_t=\psi_t \circ f \circ \phi_t^{-1}$, where $\psi_t$ and $\phi_t$ are deformations that satisfy $\psi_0=\id$ and $\phi_0=\id$. Hence, if a map is stable, then any perturbation that is near it has the same kind of singularity. In order to make precise this intuition, one has to provide a rigorous meaning of what \textit{deformation} means. This is done through the concept of \textit{unfolding}.
\begin{definition} Let $f:(\C^n,S)\rightarrow (\C^p,0)$. An $r$\textit{-parameter unfolding} of $f$ is a germ of mapping $F:(\C^n\times \C^r,S\times 0)\rightarrow (\C^p\times \C^r,0)$ of the form $F(x,t)=(\tilde{f}(x,t),t)$, and such that, if $f_t(x)=\tilde{f}(x,t)$, then $f_0=f$. 
\end{definition}
With this definition, one recovers the desired notion of \textit{deformation}, since $f_t$ are smooth maps with $f_0=0$, and that depend smoothly on $t$. When working with unfoldings, a helpful concept to consider is whether two unfoldings can be deemed equivalent. This notion is critical when distinguishing between deformations of the same map-germ.
\begin{definition} Two unfoldings $F,G:(\C^n\times \C^r,S\times 0)\rightarrow (\C^p\times \C^r,0)$ of the same germ $f:(\C^n,0)\rightarrow (\C^p,0)$ are said to be $\A$\textit{-equivalent} if there exist germs of isomorphisms $\Phi:(\C^n\times \C^r,S\times 0)\rightarrow(\C^n\times \C^r,0)$ and $\Psi:(\C^p\times \C^r,0)\rightarrow(\C^p\times \C^r,0)$ which are themselves unfoldings of the identity in $(\C^n,S)$ and $(\C^p,0)$, respectively, and such that $G=\Psi\circ F\circ \Phi^{-1}$. 
\end{definition}
Notice that the given definition of $\A$-equivalence for unfoldings is more restrictive than the general $\A$-equivalence of germs, since one asks the isomorphisms to be unfoldings of the identity themselves. If $f_t$ and $g_t$ are the deformations defined by $F$ and $G$, respectively, and $\phi_t, \psi_t$ are the deformations of the identity induced by $\Phi$ and $\Psi$, then the previous definition is equivalent to have $g_t=\psi_t\circ f_t \circ \phi_t^{-1}$. With this, one can formulate a rigorous definition for stability.
\begin{definition} An unfolding $F$ is $\A$-\textit{trivial} if it is $\A$-equivalent to the unfolding $f\times\text{id}$, which is the ``constant'' unfolding given by $(x,t)\mapsto (f(x),t)$. In these terms, one says that a map-germ $f$ is $\A$-\textit{stable} if every unfolding of $f$ is $\A$-trivial.  
\end{definition}
\begin{remark} Let us briefly comment on the notation we follow for the defined notions. As we are working under the $\A$-equivalence relation, all the terms that depend on this equivalence are prefixed by the symbol $\A$. In the next section, we will introduce another equivalence relation which has similar concepts but with different properties.
\end{remark}
The singularity theory of mappings aims to measure how far is a map-germ from being stable. The aim is to investigate which of the perturbations of $f$ are $\A$-trivial in the space of all possible perturbations. Moreover, we aim to gain a thorough understanding of the space of deformations that effectively deform the mapping. To accomplish this, we consider the space of \textit{infinitesimal deformations}, which can be defined as the space
\begin{equation*}
    ID(f)=\left\{\left.\dfrac{df_t}{dt}\right|_{t=0}: F(x,t)=(f_t(x),t) \textit{ is an unfolding of }f\right\}.
\end{equation*}
Intuitively, this captures all the possible directions in which the map-germ $f$ can be deformed. It turns out that one can identify $ID(f)$ with the vector fields defined along $f$, which are $\theta(f)=\{\xi:(\C^n,S)\rightarrow T\C^p: \xi(x)=(f(x), \tilde{\xi}(x))\}$. Indeed, associate to each element of $ID(f)$ of the form $(df_t/dt)|_{t=0}$ the vector field in $\theta(f)$ given by 
\begin{equation*}
    \xi(x)=(f(x), (df_t/dt)|_{t=0}(x)).
\end{equation*} 
Conversely, for each vector field $\xi(x)=(f(x), \tilde{\xi}(x))$, one can construct the deformation $f_t(x)=f(x)+t\tilde{\xi}(x)$, which clearly satisfies that $(df_t/dt)|_{t=0}=\tilde{\xi}$. Therefore, one can identify $ID(f)\equiv\theta(f)$. \\

Studying perturbations through vector fields offers the advantage that the space $\theta(f)$ possesses an $\O_{n,S}$-module structure, where $\O_{n,S}$ is the ring of holomorphic functions $f:(\C^n,S)\rightarrow \C$. Thus, $\theta(f)$ has a natural product defined for $f\in \O_{n,S}$ and $\xi\in \theta(f)$ given by the usual multiplication $(f\cdot \xi) (x)=f(x)\xi(x)$. This endows the space of infinitesimal deformations $\theta(f)$ with not only a complex vector space structure but also a richer module structure.\\

Our objective is to quantify the lack of stability of a germ, which entails determining how far the deformations of $\theta(f)$ are from being $\A$-trivial.
\begin{definition} Let $f:(\C^n,S)\rightarrow (\C^p,0)$. Define 
\begin{equation*}
    T\A_e f=\left\{\left. \dfrac{d(\psi_t\circ f\circ \phi_t^{-1})}{dt}\right|_{t=0} : \phi_0=\id, \psi_0=\id \right\}\subset \theta(f)
\end{equation*}
to be the \textit{extended} $\A$-\textit{tangent space of} $f$, formed by the infinitesimal deformations of $f$ that are $\A$-trivial, and let the quotient module
\begin{equation*}
    T^1_{\A_e} f=\frac{\theta(f)}{T\A_e f}
\end{equation*}
be the \textit{extended} $\A$-\textit{normal space of} $f$. Define the $\A_e$\textit{-codimension} of the germ $f$ to be
\begin{equation*}
    \codimAe(f)=\dim_{\C}T^1_{\A_e} f=\dim_{\C} \dfrac{\theta(f)}{T\A_e f}.
\end{equation*}
In this context, we say that a map-germ $f$ is $\A$\textit{-finite} if it has finite $\A_e$-codimension. 
\end{definition}
\begin{remark} The defined modules are called \textit{extended} due to the fact that, in the literature, there is a non-extended version of the $\A$-tangent space. However, the provided definition suffices for our current objectives.
\end{remark}
The following theorem establishes an important connection between the stability of a germ and its codimension:
\begin{theorem} A map-germ $f$ is stable if and only if $\codimAe(f)=0$.
\end{theorem}
\begin{proof} The direct implication follows easily, since every unfolding $F(x,t)=(f_t(x),t)$ of $f$ is $\A$-trivial, and hence $f_t=\psi_t\circ f\circ \phi_t^{-1}$ for some families of isomorphisms $\phi_t, \psi_t$ with $\phi_0=\id$, $\psi_0=\id$. Thus, $T\A_e f=\theta(f)$, forcing $\codimAe(f)=0$. The inverse implication is much more delicate, and we omit the details here. The proof can be found in Theorem 3.2 of \cite{juanjo}. 
\end{proof}
This result states that the stable germs are precisely the germs that have $T\A_e f=\theta(f)$. In general, the codimension quantifies how far is the germ $f$ to be stable. \\

Understanding the structure of the extended $\A$-tangent space is crucial for analyzing the singularity of a map-germ $f$. Thus, the following lemma is of significant importance as it offers a means of decomposing the extended $\A$-tangent space into two components.
\begin{lemma} Let $\phi_t$ and $\psi_t$ be parametrised families of isomorphisms. Then, 
\begin{equation*}
    \left.\dfrac{d(\psi_t\circ f\circ \phi_t^{-1})}{dt}\right|_{t=0}=df\circ \left(\left.\dfrac{d\phi_t^{-1}}{dt}\right|_{t=0}\right)+\left(\left. \dfrac{d\psi_t}{dt}\right|_{t=0}\right)\circ f.
\end{equation*}
\end{lemma}
The proof of the result is just a clever application of the chain rule, and can be found in Lemma 3.2 of \cite{juanjo}. 
\begin{remark}
This lemma provides a way to decompose the module $T\A_e f$ in two different parts. Indeed, if we denote $\theta_{n,S}=\{\xi:(\C^n,S)\rightarrow T\C^n: \xi(x)=(x,\tilde{\xi}(x))\}$ as the vector fields in $(\C^n,S)$, consider the map
\begin{equation*}
    tf:\theta_{n,S} \rightarrow \theta (f)
\end{equation*}
given by $tf(\xi)=df\circ \xi$, that is, $tf(\xi)(x)=(f(x),df_x\circ\tilde{\xi}\,(x))$, and the map
\begin{equation*}
    \omega f: \theta_{p,0} \rightarrow \theta (f)
\end{equation*}
defined as $\omega f(\eta)=\eta\circ f=(f, \tilde{\eta}\circ f)$. Hence, the previous lemma states that one can write 
\begin{equation*}
    T\A_e f=tf(\theta_{n,S})+\omega f(\theta_{p,0}). 
\end{equation*}
Hence, the $\A_e$-codimension of $f$ can be computed as 
\begin{equation*}
    \codimAe(f)=\dim_{\C}\dfrac{\theta(f)}{tf(\theta_{n,S})+\omega f(\theta_{p,0})}.
\end{equation*}
When performing calculations involving monogerms centered at the origin, that is, when $S={0}$, we identify the space of infinitesimal deformations $\theta_{n}$ with $\theta_{n,0}$, and we can consider it as the column space $\O_n^n$. Similarly, we identify $\theta(f)$ with $\O_n^p$ in the same way.
\end{remark}
\begin{example}\label{example:codimension} Let $f:(\C,0)\rightarrow (\C^2,0)$ be the map $f(x)=(x^2,x^3)$. Let us compute the extended $\A$-tangent space of $f$. Firstly, we shall compute the space $\omega f(\theta_2)$, which are elements of the form 
\begin{equation*}
    \begin{bmatrix}
        a(x^2,x^3) \\ b(x^2,x^3)
    \end{bmatrix},
\end{equation*}
where $a,b\in \O_2$. Notice that any monomial $x^k$ except for $x$ can be expressed as a function of $x^2$ and $x^3$. Therefore, 
\begin{equation*}
    \omega f(\theta_2)+\sp_{\C}\Biggl\{ \begin{bmatrix}x \\ 0 \end{bmatrix}, \begin{bmatrix}0 \\ x \end{bmatrix} \Biggr\}=\theta (f).
\end{equation*} 
Moreover, notice that 
\begin{equation*}
    tf\left(\dfrac{\partial}{\partial x}\right)=\begin{bmatrix} 2x \\ 3x^2  \end{bmatrix},
\end{equation*}
and thus 
\begin{equation*}
    \begin{bmatrix}x \\ 0 \end{bmatrix}=\dfrac{1}{2} \begin{bmatrix}2x \\ 3x^2 \end{bmatrix}-\dfrac{3}{2}\begin{bmatrix}0 \\ x^2 \end{bmatrix}\in tf(\theta_1)+\omega f(\theta_2)=T\A_ef.
\end{equation*}
Hence, one has that
\begin{equation*}
    T\A_ef +\sp_{\C} \Biggl\{\begin{bmatrix}0 \\ x \end{bmatrix} \Biggr\}=\theta (f),
\end{equation*}
and the term $\begin{bmatrix} 0 \\ x\end{bmatrix}\notin T\A_e f$ since the entries in the second row of $T\A_e f$ should have order $\geq 2$. It follows that 
\begin{equation*}
    \codimAe (f)=\dim_{\C} \dfrac{\theta(f)}{T\A_e f}=\dim_{\C} \sp_{\C} \Biggl\{\begin{bmatrix}0 \\ x \end{bmatrix} \Biggr\}=1.
\end{equation*}
\end{example}

The next crucial definition is the concept of finite determinacy of mappings. Recall that the $k$-jet of $f$ at a point $x$ is the Taylor polynomial mapping of order $k$ centered at $x$, and it is denoted as $j^k f(x)$. 
\begin{definition}
    For $k\in \N$, we say that $f$ is $k$\textit{-determined for} $\A$\textit{-equivalence} if, for every map $g$ such that the $k$-jets at the origin $j^kf(0)$ and $j^kg(0)$ are $\A$-equivalent, then $f$ and $g$ are $\A$-equivalent. We say that $f$ is \textit{finitely determined for} $\A$\textit{-equivalence} if it is $k$-determined for some $k\in \N$.
\end{definition}
Let us denote $\m_{n,S} =\{f\in \O_{n,S}: f(S)=0\}$ as the function-germs that vanish at the points of $S$. An important result from Mather states that finite determinacy is equivalent to $\A$-finiteness. 
\begin{theorem}(Mather's finite determinacy Theorem) For a map-germ $f:(\C^n,S)\rightarrow (\C^p,0)$, the following are equivalent:
    \begin{enumerate}
        \item $f$ is $\A$-finite.
        \item $f$ is finitely determined for $\A$-equivalence. 
        \item $\m_{n,S}^k \theta(f)\subset T\A_e f$ for some $k\in \N$. 
    \end{enumerate}
\end{theorem}
In fact, there is a close relationship between the number $k$ that appears in the theorem and the degree of determinacy of the function. Since this detail and the proof of the result will not be relevant for our purposes, we refer the reader to theorem 6.2 in the reference \cite{juanjo}.\\

Recall that, for a mapping $f:U\subset \C^n \rightarrow \C^p$, we define the set of critical points as $C(f)=\{x\in U: df_x \text{ is not surjective}\}$, and the discriminant of $f$ as $\Delta(f)=f(C(f))$. A crucial result for the $\A$-equivalence relation is the celebrated Mather-Gaffney criterion, which states that a map-germ is $\A$-finite if and only if it has isolated instability. 
\begin{theorem} (Mather-Gaffney) A map-germ $f:(\C^n,S)\rightarrow (\C^p,0)$ is $\A$-finite if and only if there is a small enough representative $f:U\rightarrow V$ such that $f^{-1}(0)\cap C(f)=S$ and that $f:U\setminus f^{-1}(0)\rightarrow V\setminus\{0\}$ is locally stable. 
\end{theorem}
For a proof of the result, see theorem 4.5 of \cite{juanjo}. \\

In what follows, let us give a condition for an unfolding to capture all possible perturbations that arise near a singularity. 
\begin{definition} Let $f:(\C^n,S)\rightarrow (\C^p,0)$ be a map-germ, and $F:(\C^n\times\C^r,S\times 0)\rightarrow (\C^p\times\C^r,0)$ be an unfolding of the form $F(x,t)=(f_t(x),t)$, where $f_0=f$. Given a map-germ $h:(\C^d,0)\rightarrow (\C^r,0)$, define the \textit{pull-back} of $F$ by $h$ as the unfolding $h^*F:(\C^n\times\C^d,S\times 0)\rightarrow (\C^p\times\C^d,0)$ defined by $(x,u)\rightarrow (f_{h(u)}(x), u)$.
\end{definition}
The intuition behind this definition is that $h^*F$ is an unfolding where the information of $F$ is deformed perturbing the parameter space. 
\begin{definition} An unfolding $F$ of $f$ is $\A$-\textit{versal} if any other unfolding of $f$ is $\A$-equivalent to a pull-back $h^*F$ for some base-change $h$. We say that $F$ is $\A$-\textit{miniversal} if it is an $\A$-versal unfolding with the property that there are no $\A$-versal unfoldings with less parameters than $F$. 
\end{definition}
Intuitively, an  $\A$-versal unfolding contains all other possible unfoldings through perturbations of the parameters. Furthermore, an $\A$-miniversal one has the smallest possible number of parameters. \\

Although the given definition has a strong underlying intuition, it is extremely inoperative in practice. Indeed, in order to check if an unfolding $F$ is $\A$-versal, one is expected to show that any other unfolding is $\A$-equivalent to a pull-back $h^*F$. The following theorem gives a clear method to compute an $\A$-miniversal unfolding of a map-germ.
\begin{theorem}\label{theorem:versality} Let $f:(\C^n,S)\rightarrow (\C^p,0)$ be a map-germ. An unfolding $F:(\C^n\times \C^r,S\times 0)\rightarrow (\C^p\times\C^r,0)$ of the form 
    $$F(x,u)=(f(x)+u_1g_1(x)+\ldots+u_rg_r(x),u)$$
is $\A$-versal if and only if 
\begin{equation*}
    T\A_e f+ \sp_\C \{g_1, \ldots, g_r\}=\theta(f). 
\end{equation*}
\end{theorem}
The proof of the theorem can be found in \cite{juanjo}, in theorem 5.1. An important corollary of this result is 5.1 of \cite{juanjo}, which states the following:
\begin{corollary}\label{corollary:miniversal} Let $f:(\C^n,S)\rightarrow (\C^p,0)$ be a map-germ. Then,
\begin{enumerate}
    \item $f$ admits an $\A$-versal unfolding if and only if $f$ is $\A$-finite. 
    \item $\codimAe(f)$ is the number of parameters in an $\A$-miniversal unfolding of $f$. 
\end{enumerate}
\end{corollary}
\begin{example} Let us consider the map-germ $f(x)=(x^2,x^3)$ of Example \ref{example:codimension}. Recall that
\begin{equation*}
    T\A_ef +\sp_{\C} \Biggl\{\begin{bmatrix}0 \\ x \end{bmatrix} \Biggr\}=\theta (f),
\end{equation*}
and hence the unfolding $F(x,t)=(x^2,x^3+tx,t)$ is $\A$-versal by Theorem \ref{theorem:versality}. Since $\codimAe(f)=1$ is the number of parameters of $F$, it follows by Corollary \ref{corollary:miniversal} that $F$ is $\A$-miniversal.
\end{example}
This completes a concise overview of the basic terminology used to classify map-germs under $\A$-equivalence and determine its possible deformations. However, the extended $\A$-tangent space $T\A_e f$, which is the central algebraic object of this equivalence, has a complex structure that can be difficult to fully understand in general. In the following section, we explore an alternative equivalence that offers a solution to this issue. Although the definition of this new equivalence may seem unclear at first, its associated tangent space possesses a more straightforward structure that facilitates computation and comprehension. This new relation not only aids in understanding map-germs through a different perspective, but also proves immensely useful in the study of the $\A$-equivalence classification of map-germs. 
\subsection{$\K$-equivalence of mappings}
There exist several reasons why $\A$-equivalence is insufficient for properly comprehending singularities. As we have commented on in the last section, the extended tangent space associated with $\A$-equivalence can be very challenging to directly calculate for a general map-germ. To overcome this issue, an analogous tangent space for $\K$-equivalence has been developed and serves as a useful intermediate tool for determining the $\A$-tangent space. Furthermore, when analyzing singularities of analytic space germs, $\K$-equivalence arises naturally. Specifically, we will show in the next chapter that the isomorphism class of certain analytic germs, namely the isolated complete intersection singularities, can be entirely determined by the $\K$-equivalence class of the mappings defining them. As a result, $\K$-equivalence is an indispensable prerequisite for this project. This section will explore this distinct equivalence for mappings, highlighting both the similarities and differences with $\A$-equivalence.
\begin{definition} Let $f,g:(\C^n,S)\rightarrow (\C^p,0)$ be germs of mappings. We say that $f$ and $g$ are $\K$\textit{-equivalent} if there exists an isomorphism $\Psi:(\C^n\times\C^p,S\times 0)\rightarrow (\C^n\times\C^p,S\times 0)$ of the form $\Psi(x,y)=(\phi(x),\psi(x,y))$ such that $\psi(x,0)=0$ for all $x$, and 
    \begin{equation*}
        g(\phi(x))=\psi(x,f(x)). 
    \end{equation*}
\end{definition}
\begin{remark} At first sight, this definition could look artificial and technical. In order to provide a clearer intuition of the concept, notice that, if $f$ and $g$ are $\K$-equivalent, then 
\begin{equation*}
    (\phi(x), g(\phi(x)))=\Psi(x,f(x)).
\end{equation*}
Hence, $\Psi$ is an isomorphism that sends the graph of $f$ to the graph of $g$. In other words, if $f$ and $g$ are $\K$-equivalent, there exists an isomorphism $\Psi$ in $(\C^n\times\C^p,S\times 0)$ of the previous form that carries the graph of $f$ to the graph of $g$, and leaves $(\C^n\times \{0\},S\times 0)$ invariant. \\
Furthermore, let us consider the \textit{graph embedding} of a function $f:(\C^n,S)\rightarrow (\C^p,0)$ to be the map $\gr(f):(\C^n,S)\rightarrow (\C^n\times \C^p,S\times 0)$ given by $\gr(f)(x)=(x,f(x))$. Then, the maps $f,g$ are $\K$-equivalent if and only if $\gr(g)\circ \phi = \Psi \circ \gr(f)$. Therefore, $f$ and $g$ are $\K$-equivalent if and only if the diagram 
    \begin{center}
        \begin{tikzcd}
          (\C^n,S) \arrow[r, "\gr(f)"] \arrow[d, "\phi"]
          & (\C^n\times\C^p,S\times 0) \arrow[d, "\Psi"] \\
          (\C^n,S) \arrow[r, "\gr(g)"]
          & (\C^n\times\C^p,S\times 0)
          \end{tikzcd}
      \end{center}
      commutes. This means that the functions $f$ and $g$ are $\K$-equivalent if and only if their graph embeddings are $\A$-equivalent through the pair $(\phi,\Psi)$. Thus, the study of the map-germs $f:(\C^n,S)\rightarrow (\C^p,0)$ under $\K$-equivalence is strongly linked with the $\A$-equivalence relation, but it rather has a more geometrical flavour.  
\end{remark}

\begin{remark} In a similar way as it was done in remark \ref{remark:Aequivalence}, $\K$-equivalence of mappings can be understood through an action of a group of isomorphisms. Indeed, consider the group $\K$ of the isomorphisms $\Psi:(\C^n\times \C^p,S\times 0)\rightarrow (\C^n\times\C^p,S\times 0)$ that have the form $\Psi(x,y)=(\phi(x), \psi(x,y))$, and $\psi(x,0)=0$ for all $x$. The action of $\K$ over the space of mappings $f:(\C^n,S)\rightarrow (\C^p,0)$ is given by 
\begin{equation*}
    (\Psi\cdot f) (x)=\psi(\phi^{-1}(x), f\circ \phi^{-1}(x)). 
\end{equation*}
Hence, it follows that $f$ and $g$ are $\K$-equivalent if and only if $g=\Psi\cdot f$ for some $\Psi\in \K$. 
\end{remark}
It seems natural that many definitions stated for the $\A$-equivalence have a parallel version in this new setting. In what follows, we give a brief outline of how can the notions of the first subsection be adapted to the $\K$-equivalence of mappings. 
\begin{definition}\textcolor{white}{.}
\begin{enumerate}
    \item Two unfoldings $F,G:(\C^n\times \C^r,S\times 0)\rightarrow (\C^p\times \C^r,0)$ of the same germ $f:(\C^n,S)\rightarrow (\C^p,0)$ are said to be $\K$\textit{-equivalent} if there exists a germ of isomorphism $\Psi:(\C^n\times \C^p \times \C^r,S\times 0)\rightarrow(\C^n\times \C^p\times \C^r,S\times 0)$ of the form $\Psi(x,y,t)=(\phi_t(x),\psi_t(x,y),t)$, where $\psi_t(x,0)=0$ for all $(x,t)$, which is an unfolding of the identity in $(\C^n\times\C^p,S\times 0)$, and such that $G=\Psi\cdot F$. In other words, if $f_t$ and $g_t$ are the deformations associated to the unfoldings $F$ and $G$, respectively, one has that
    \begin{equation*}
        (g_t\circ \phi_t) (x)=\psi_t(x,f_t(x)). 
    \end{equation*}
    \item  An unfolding $F$ is $\K$-\textit{trivial} if it is $\K$-equivalent to the unfolding $f\times\text{id}$. In these terms, one says that a map-germ $f$ is $\K$-\textit{stable} if every unfolding of $f$ is $\K$-trivial.   
\end{enumerate}
\end{definition}
As it was done for $\A$-equivalence, one can try to measure how far are infinitesimal deformations of a map-germ $f$ to be $\K$-trivial. In particular, we have the following concepts:
\begin{definition} Let $f:(\C^n,S)\rightarrow (\C^p,0)$. Define the module 
    \begin{equation*}
        T\K_e f=tf(\theta_n)+f^*(\m_p)\theta(f)\subset \theta(f)
    \end{equation*}
    as the \textit{extended} $\K$-\textit{tangent space} of $f$, where $f^*(\m_p)$ is the ideal in $\O_n$ generated by the functions $f_1,\ldots, f_p$. Let the quotient module 
    \begin{equation*}
        T^1_{\K_e} f=\frac{\theta(f)}{T\K_e f}
    \end{equation*}
    be the $\K$-\textit{normal space of} $f$. Define the $\K_e$\textit{-codimension} of the germ $f$ to be
    \begin{equation*}
        \codimKe(f)=\dim_{\C}T^1_{\K_e} f=\dim_{\C} \dfrac{\theta(f)}{T\K_e f}.
    \end{equation*}
    In this context, we say that a map-germ $f$ is $\K$\textit{-finite} if it has finite $\K_e$-codimension.
\end{definition}
In the same way we did for the $\A$-tangent space, one can describe heuristically the $\K$-tangent space of a function $f$ as the directions of $\K$-trivial deformations of $f$. However, we omit this detail and refer it to Lemma 4.1 of \cite{juanjo}. 
\begin{remark} In general, the obtention of the $\A$-tangent space of a map-germ is a complicated task. In this sense, the module $T\K_e f$ plays an auxiliary role for determining $T\A_e f$, since the $\K$-tangent space is easier to obtain, and its structure plays a key role for determining $T\A_e f$. For more details on this, see Section 4.2 of \cite{juanjo}. 
\end{remark}
Notice that the same definition of finite determinacy can be given in this new setting:
\begin{definition}
    For $k\in \N$, we say that $f$ is $k$\textit{-determined for} $\K$\textit{-equivalence} if, for every map $g$ such that the $k$-jets at the origin $j^kf(0)$ and $j^kg(0)$ are $\K$-equivalent, then $f$ and $g$ are $\K$-equivalent. We say that $f$ is \textit{finitely determined for} $\K$\textit{-equivalence} if it is $k$-determined for some $k\in \N$.
\end{definition}
The same result for finite determinacy that holds in $\A$-equivalence turns out to hold in general for all Mather's groups. In particular, the result also holds for $\K$-equivalence:
\begin{theorem}(Mather's finite determinacy Theorem) For a map-germ $f:(\C^n,S)\rightarrow (\C^p,0)$, the following are equivalent:
    \begin{enumerate}
        \item $f$ is $\K$-finite.
        \item $f$ is finitely determined for $\K$-equivalence. 
        \item $\m_{n,S}^k \theta(f)\subset T\K_e f$ for some $k\in \N$. 
    \end{enumerate}
\end{theorem}
Both results for $\A$ and $\K$-equivalence are a particular case of Theorem 6.2 of \cite{juanjo}. \\

The forthcoming statement presents an analogous criterion to the Mather-Gaffney criterion for $\A$-equivalence. However, for $\K$-equivalence, the criterion has a more geometric flavor, as will be shown in the connection between the $\K$-equivalence of map-germs and the isomorphism class of isolated complete intersection singularities. This is why the following result is commonly known as the ``geometric criterion'':
\begin{theorem} [Geometric criterion] A mapping $f:(\C^n,S)\rightarrow (\C^p,0)$ is $\K$-finite if and only if it is finite-to-one on its critical set. 
\end{theorem}
Regarding versality of unfoldings, the definitions given for $\A$-equivalence have an analogous version:
\begin{definition} An unfolding $F$ of $f$ is $\K$-\textit{versal} if any other unfolding of $f$ is $\K$-equivalent to a pull-back $h^*F$ for some base-change $h$. We say that $F$ is $\K$-\textit{miniversal} if it is a $\K$-versal unfolding with the property that there are no $\K$-versal unfoldings with less parameters than $F$. 
\end{definition}
The practical results to find a $\K$-miniversal unfolding for a map-germ hold \textit{mutatis mutandis} as for the case of $\A$-equivalence. 
\begin{theorem}\label{theorem:Kversality} Let $f:(\C^n,S)\rightarrow (\C^p,0)$ be a map-germ. An unfolding $F:(\C^n\times \C^r,S\times 0)\rightarrow (\C^p\times\C^r,0)$ of the form $F(x,u)=(f(x)+u_1g_1(x)+\ldots+u_rg_r(x),u)$ is $\K$-versal if and only if 
    \begin{equation*}
        T\K_e f+ \sp_\C \{g_1, \ldots, g_r\}=\theta(f). 
    \end{equation*}
    \end{theorem}
    \begin{corollary}\label{corollary:Kminiversal} Let $f:(\C^n,S)\rightarrow (\C^p,0)$ be a map-germ. Then,
    \begin{enumerate}
        \item $f$ admits a $\K$-versal unfolding if and only if $f$ is $\K$-finite. 
        \item $\codimKe(f)$ is the number of parameters in a $\K$-miniversal unfolding of $f$. 
    \end{enumerate}
    \end{corollary}
\begin{example} Let us compute a $\K$-miniversal unfolding of the map-germ $f:(\C^2,0)\rightarrow (\C^2,0)$ given by $f(x,y)=(x^2,y^2)$. First of all, notice that $tf(\theta_2)$ is the submodule of $\theta(f)$ generated by the vectors 
\begin{equation*}
    \begin{bmatrix}        2x \\ 0    \end{bmatrix},    \begin{bmatrix}        0 \\ 2y    \end{bmatrix},
\end{equation*}
and $f^*(\m_2)\theta(f)=(x^2,y^2)\,\theta(f)$. Therefore,
\begin{equation*}
    T\K_e f+\sp_{\C}\left\{ \begin{bmatrix}        y \\ 0    \end{bmatrix},    \begin{bmatrix}        0 \\ x    \end{bmatrix} \right\}=\theta(f).
\end{equation*}
Moreover, $\codimKe (f)=2$, since the previous two vectors do not lie in $T\K_e f$. Hence, $F(x,y,u,v)=(x^2+uy,y^2+vx,u,v)$ is a $\K$-miniversal unfolding of $f$. 
\end{example}
A particular fact for $\K$-equivalence of mappings is that it has a complete algebraic invariant, namely the \textit{semi-local algebra} of a map-germ $f$. Recall that an \textit{algebra} over $\C$ is a $\C$-vector space $A$ equipped with a ring structure with identity, in such a way that its product operation as a ring is compatible with its scalar product as a vector space. In other words, for all $\lambda \in \C$ and all $x,y\in A$ one has that 
\begin{equation*}
    (\lambda x)y=\lambda (xy)=x(\lambda y).
\end{equation*}
In this context, we say that two algebras $A$ and $B$ are isomorphic if there exists a map $\phi:A\rightarrow B$ that is both an isomorphism of vector spaces and an isomorphism rings. \\

In particular, we are interested in \textit{analytic} algebras, which are $\C$-algebras of the form $A=\O_{n,S}/I$ for some ideal $I\subset \O_{n,S}$. 

\begin{definition}
    Let $f:(\C^n,S)\rightarrow (\C^p,0)$ be a map-germ. Define the \textit{semi-local algebra} of $f$ as the analytic algebra 
    \begin{equation*}
        Q(f)=\dfrac{\O_{n,S}}{f^*(\m_p)\O_{n,S}}=\dfrac{\O_{n,S}}{(f_1,\ldots, f_p)},
    \end{equation*}
    where $(f_1,\ldots,f_p)$ stands for the ideal in $\O_{n,S}$ generated by the coordinate functions of $f$. 
\end{definition}
This algebra turns out to determine completely the $\K$-equivalence class of a mapping, as the following theorem states:
\begin{theorem} The map-germs $f,g:(\C^n,S)\rightarrow (\C^p,0)$ are $\K$-equivalent if and only if their semi-local algebras $Q(f)$ and $Q(g)$ are isomorphic. 
\end{theorem}
\begin{remark} In the following section, we will show that the algebra $Q(f)$ is, in fact, the algebra of germs in the analytic space given by $(f^{-1}(0),S)$. Hence, the classification of mappings under $\K$-equivalence is tightly linked with the classification of some analytic spaces. 
\end{remark}

\newpage
\thispagestyle{empty}
\textcolor{white}{blank page}
\newpage
\section{Isolated complete intersection singularities (\textsc{icis})}
One of the primary focuses of this project is to study a particular type of analytic spaces with exceptional properties known as isolated complete intersection singularities (\icis\!\!). These spaces possess two distinct characteristics that set them apart from other analytic spaces. Firstly, the number of equations required to define an \icis is equal to its codimension in the ambient space, which is why they are referred to as complete intersections. Secondly, the singular locus of an \icis is the simplest nontrivial possible case, consisting of a single point, which is why they are referred to as isolated singularities. The exceptional richness of the analytic spaces where both conditions are met has led to their in-depth study by singularists in the last years. In this chapter, we provide an overview of their fundamental properties and their possible deformations. We also introduce two crucial numerical invariants, the Milnor and Tjurina numbers, to establish a profound relationship between them.

\subsection{Basic definitions and examples}
This first section aims to provide an overview and introduce the fundamental concepts related to isolated complete intersection singularities. We presume that the reader has a basic understanding of algebraic and analytic geometry. Specifically, we assume that the reader is familiar with the concept of the dimension of an analytic space $(\C^n,S)$. Additionally, we presume some background knowledge regarding commutative algebra's elementary properties.
\begin{definition} Let $(X,S)\subset (\C^n,S)$ be a germ of complex space of dimension $d$. We say that $(X,S)$ is \textit{smooth} or that $x_0$ is a \textit{regular point} of $X$ whenever $(X,S)$ is isomorphic to $(\C^d,S)$. In case $(X,x_0)$ is not regular, where $x_0\in S$, we say that $x_0$ is a singular point of $X$, or that $(X,x_0)$ is singular. 
\end{definition}
Fix some representative $X\subset \C^n$ and let $g_1, \ldots, g_s$ be equations defined on $U\subset \C^n$ such that $X\cap U=V(g_1, \ldots, g_s)$. We can therefore consider $\Sigma$ as the set of points $x\in X\cap U$ such that $(X,x)$ is singular. The set-germ $(\Sigma, S)$ is called the \textit{singular locus} of $(X,S)$, and it is denoted as Sing$\,(X,S)$. In order to determine this set-germ in practice, the following theorem becomes crucial:

\begin{theorem}[Jacobian criterion] Assume that $\dim (X,x)=n-c$ for each $x\in X$. Then, Sing$\,(X,S)$ is the subset-germ of $(X,S)$ defined as the vanishing set of the $c$-minors of the jacobian matrix $dg$ of $g=(g_1, \ldots, g_s)$. 
\end{theorem}
In particular, the previous theorem asserts that Sing$\,(X,S)$ is an analytic set, and provides the neccesary equations to define it. Its proof is rather nontrivial, and it is reffered to more advanced references. \\

Studying analytic sets without any restriction on their singular locus can be a challenging task to tackle. Therefore, to simplify the situation, researchers typically focus on studying analytic sets with \textit{isolated singularities}, where the singular locus is either a point or the empty set. Such sets are considered to be the easiest examples of analytic sets in terms of their singular locus. Many of the results presented in the literature are limited to complex spaces with isolated singularities. The theory of these spaces is much more extensive and richer than that of general analytic sets.
\begin{example}\label{example:IS_NCI} Let $(X,0)\subset (\C^4,0)$ be the germ defined by the equations $xz-y^2, yw-z^2, xw-yz$. With \singular\!, one can easily check that $\dim (X,0)=2$. Hence, the singular locus is formed by the points of $(X,0)$ where the $2$-minors of the jacobian matrix 
\begin{equation*}
    \begin{pmatrix}
        z & -2y & x & 0 \\
        0 & w & -2z & y \\
        w & -z & -y & x 
    \end{pmatrix}
\end{equation*}
vanish. After some calculations, it easily follows that Sing$\,(X,0)=({0},0)$. Therefore, $(X,0)$ is an isolated singularity.     
\end{example}
\begin{definition}
    Let $(X,S)\subset (\C^n,S)$ be a germ of complex space defined by an ideal $I\subset \O_{n,S}$, and let $k$ be the minimum number of generators of $I$. We say that $(X,S)$ is a \textit{complete intersection} if $\dim (X,S)=n-k$. 
\end{definition}
\begin{example} Any hypersurface is a complete intersection, since they are $n-1$ dimensional spaces defined by a single equation. In general, the number of generators of an ideal $I$ is greater than or equal to the codimension of the associated analytic space. However, there exist examples where the minimum number of generators of $I$ is strictly greater than the codimension of the space. For instance, the set-germ of example \ref{example:IS_NCI} is not a complete intersection, since it is a two-dimensional germ that cannot be generated by two equations in $\C^4$. Another famous example is the curve defined by the image of  $t\mapsto (t^3,t^4,t^5)$ in $(\C^3,0)$. This curve is generated by the three equations $x^4-y^3$, $x^5-z^3$, and $y^5-z^4$. It can be shown that this set cannot be described by only two equations. Therefore, this curve is not a complete intersection.
\end{example}
\begin{definition} A complex space that is a complete intersection and that has isolated singularity is called an \textit{isolated complete intersection singularity}, and it is shortened as \icis\!. 
\end{definition}
\begin{example} The example given in \ref{example:IS_NCI} is not an \icis due to the fact that it is not a complete intersection. However, it can be shown that the germ $(X,0)\subset (\C^3,0)$ defined by the equations $x^2-y^3, x^3-z^2$ is a complete intersection (just by checking that its dimension is 1) and has an isolated singularity. Hence, $(X,0)$ is an \icis \!.
\end{example}
As mentioned earlier, the purpose of these definitions is to narrow our focus to certain types of spaces that exhibit particular forms of well-behaved behavior. Specifically, we are interested in studying spaces with a singular locus that is as simple as possible, yet nontrivial (isolated singularities), and whose defining equations are minimal, i.e., equal to its codimension (complete intersections). \\

Nonetheless, many other types of analytic spaces with special properties have been defined along the literature. In this project, a key concept in this regard is the Cohen-Macaulay property, which is typically formulated in its general version for modules. 
\begin{definition} Let $(R,\m)$ be a local ring and $M$ be an $R$-module. 
\begin{enumerate}
    \item A sequence $f_1, \ldots, f_r$ of elements in $\m$ is called a \textit{regular sequence} of $M$ if $f_1$ is not a zerodivisor of $M$, and $f_i$ is not a zerodivisor of $M/(f_1, \ldots, f_{i-1})M$ for $i=2, \ldots, r$. 
    \item Let $I\subset IR$ be an ideal with $IM\neq M$. The $I$\textit{-depth} of $M$, $\depth (I,M)$ is the maximal length of a regular sequence of $M$ in $I$. If $IM=M$ we set $\depth (I,M)=\infty$. 
    \item The \textit{depth of}$M$, $\depth M$ is the maximal length of a regular sequence of $M$, that is, $\depth M =\depth (\m, M)$. 
\end{enumerate}
If the ring $R$ is Noetherian, it can be shown that $\depth M\leq \dim M$ in any case. When equality holds, we say that the module $M$ is \textit{Cohen-Macaulay}. A ring $R$ is called \textit{Cohen-Macaulay} whenever $R$ is Cohen-Macaulay as an $R$-module. Furthermore, if an analytic space $(X,x_0)\subset (\C^n,x_0)$ is defined by the ideal $I\subset \O_{n,x_0}$, we say that it is \textit{Cohen-Macaulay} whenever $\O_{X,x_0}=\O_{n,x_0}/I$ is Cohen-Macaulay as a ring. In the case of multigerms, we say that $(X,S)$ is \textit{Cohen-Macaulay} if $(X,x_0)$ is Cohen-Macaulay for any $x_0\in S$.  
\end{definition}
\begin{remark} The definition of a Cohen-Macaulay module is a crucial concept in the heart of the project, which is chapter 5. Nevertheless, in this chapter, we center our attention primarily on the definition of Cohen-Macaulay analytic spaces, and its relation with isolated complete intersection singularities. 
\end{remark}
The definition of the Cohen-Macaulay property may appear less straightforward than that of complete intersection singularities. In order to provide examples of this concept, let us present some results concerning this property.
\begin{theorem} Any \icis is Cohen-Macaulay. 
\end{theorem}
This provides a great amount of examples of Cohen-Macaulay spaces. In order to determine some example of an analytic space that is not Cohen-Macaulay, let us state another relevant result.
\begin{theorem} Any Cohen-Macaulay analytic space is equidimensional. That is, all its irreducible components have the same dimension.
\end{theorem}
This implies that, for example, the space of $(\C^3,0)$ that is the union of the line $x=y=0$ and the plane $z=0$ is not Cohen-Macaulay. However, many examples of equidimensional spaces that are not Cohen-Macaulay can be found. 

\subsection{Deformations of \icis}
In this section, we introduce the key ideas related to the deformations of an \icis\!. These concepts are presented in a manner similar to that in the first chapter, emphasizing a profound parallelism with the deformations in the $\K$-equivalence relation. Building on the definitions outlined in the previous section, we describe an \icis $(X,S)\subset (\C^N,S)$ of dimension $n$ as the fiber of a mapping.
\begin{equation*}
    g=(g_1, \ldots, g_{N-n}):(\C^N,S)\rightarrow (\C^{N-n},0),
\end{equation*}
\textit{i.e.}, $(X,S)=V(g_1, \ldots, g_{N-n})$. Moreover, the function $g$ is necessarily $\K$-finite by the geometric criterion. Conversely, the fiber of any $\K$-finite mapping of the form stated above is an \icis of dimension $n$. Hence, there exists a strong parallelism between \icis and $\K$-finite mappings.

It is worth noting that $\O_{X,S}$, the algebra of regular functions of $(X,S)$, can be expressed as $\O_{N,S}/(f_1,\ldots,f_{N-n})=Q(f)$. That is, the algebraic invariant of $(X,S)$ coincides with the local algebra of $f$. Interestingly, this algebraic invariant is complete in both settings, meaning that it determines $(X,S)$ up to isomorphism and $g$ up to $\K$-equivalence. As a result, the parallelism goes beyond mere similarities between the concepts.
\begin{proposition}
    Let us denote
    \begin{align*}
        \mathfrak{M}&=\{\K\text{-classes of finite mappings }(\C^N,S)\rightarrow (\C^{N-n},0)\}, \\
        \mathfrak{I}&=\{ \text{isomorphism classes of \icis in }(\C^N,S)\text{ of dimension }n \}.
    \end{align*}
     The correspondence $\mathfrak{M}\rightarrow \mathfrak{I}$ that sends the class of each mapping $g$ to the class of the set-germ $(g^{-1}(0),S)$ is a bijection. 
\end{proposition}
In the subsequent discussion, we present the theory of deformations of \icis\!. To facilitate this, we first introduce an important algebraic object that provides a measure of the extent to which a space can be deformed. For an \icis $(X,S)$ which is the fiber of a $\K$-finite mapping $g$, we define
\begin{equation*}
    T^1_{X,S} = \dfrac{\O_{X,S}^{N-n}}{JM_{X}(g)}
\end{equation*}
to be the \textit{module of deformations of} $(X,S)$, where $JM_X(g)$ is the submodule of $\O_{X,S}^{N-n}$ generated by the classes of $\partial g / \partial x_i$ for $i\in \{1,\ldots, n\}$. It is therefore immediate that 
\begin{equation*}
    T^1_{X,S}\cong \dfrac{\theta(g)}{T\K_e g}=T^1_{\K_e}g.
\end{equation*}
In other words, the module of deformations of an \icis is isomorphic to the extended $\K$-tangent module of its defining equation. We define $\tau (X,S)=\dim_\C T^1_{X,S}$ to be the \textit{Tjurina number} of $(X,S)$, which is equal to $\K_e$-codimension of the mapping $g$. It can be shown that this number is invariant under isomorphism. 
\begin{remark} It is worth emphasizing that, if $(X,S)$ is a hypersurface (that is, it is given by a single equation $g:(\C^N,S)\rightarrow (\C,0)$, then the Tjurina number can be easily computed as
\begin{equation*}
    \tau (X,S)=\dim_\C \dfrac{\O_{N}}{(g)+J(g)},
\end{equation*}
since $T\K_e g =(g)+J(g)$ in this case. 
\end{remark}
In order to link the previous algebraic object with a more geometrical notion of deformation, let us introduce the following definition
\begin{definition} An $r$\textit{-parameter deformation} of an \icis $(X,S)$ is a space $(\X,T)$ and a mapping $\pi:(\X,T)\rightarrow (\C^r,0)$, where $\dim (\X,T) = \dim (X, S)+r$ and where $(X,S)$ is isomorphic to the fibre $\pi^{-1}(0)$.
\end{definition}
\begin{remark} It is important to note that the definition of $(X,S)$ is given up to isomorphism. Thus, one may identify $(X,S)$ directly with $\pi^{-1}(0)$, assuming that $(X,S)$ is contained within $(\X,S)$ and $S=T$. In these terms, an $r$-parameter deformation of $(X,S)$ is a space $(\X,S)$ containing $(X,S)$ with $\dim (\X,S)=\dim (X,S)+r$, equipped with a projection mapping $\pi:(\X,S)\rightarrow (P,0)$, such that $X=\pi^{-1}(0)$. The term "deformation" derives from the idea that the fibers $X_t=\pi^{-1}(t)$ correspond to the natural concept of "deformation" of $X_0=X$, and hence the space $\X=\bigcup_{t}X_t$ contains the deformations of $X$.
\end{remark}
With this, we define a \textit{morphism of deformations} $\pi:(\X,S)\rightarrow (\C^r,0)$ and $\pi':(\X',S)\rightarrow (\C^s,0)$ to be morphisms $g:(\X,S)\rightarrow (\X',S)$ and $h:(\C^r,0)\rightarrow (\C^s,0)$ that satisfy that $g|_{X}=g'|_{X}$ and that the diagram
    \begin{center}
        \begin{tikzcd}
          (\X,S) \arrow[r, "\pi"] \arrow[d, "g"]
          & (\C^r,0) \arrow[d, "h"] \\
          (\X',S) \arrow[r, "\pi'"]
          & (\C^s,0)
          \end{tikzcd}
      \end{center}
commutes. \\

Rewrite the following text with an academic natural english: With this, we say that a deformation $f$ is \textit{versal} if any other deformation $f'$ admits a morphism from $f'$ to $f$. If there exists a unique matrix $dh_0$ among all possible morphisms, then we call $f$ \textit{miniversal}. To compute versal and miniversal deformations of an \icis, we present the following theorem, which establishes a relationship between the versal deformations of \icis and the versal unfoldings under $\K$-equivalence of the mappings that define them.

\begin{theorem} If $(X,S)\subset (\C^N,S)$ is an \icis defined by the $\K$-finite mapping $f:(\C^N,S)\rightarrow (\C^{N-n},0)$, and $F(x,u)=(f_u(x),u)$ is an $r$-parameter unfolding of $f$, the following are equivalent:
\begin{enumerate}
    \item $F$ is a versal deformation of $(X,S)$, where $(X,S)$ is identified with its image by $x\mapsto (x,0)$ in $\C^N\times \C^r$. 
    \item $F$ is $\A$-stable.
    \item $\tilde{F}(x,u,v)=(f_u(x)+v,u,v)$ is a $\K$-versal unfolding of $f$.
\end{enumerate}
Furthermore, $\K$-miniversal unfoldings of $f$ correspond to miniversal deformations of $(X,S)$.
\end{theorem}
Hence, the comprehension of the deformations of an \icis is completely equivalent to the understanding of the unfoldings of its associated mapping with respect to $\K$-equivalence. 
\begin{remark} In particular, this theorem implies that $\tau(X,S)$ is the number of parameters of a miniversal deformation of $(X,S)$. 
\end{remark}
\begin{example} Let $(X,0)\subset (\C^3,0)$ be the set-germ defined by the mapping $f:(\C^3,0)\rightarrow (\C^2,0)$ given by $f(x,y,z)=(x^2+y^2+z^2,xy)$. A straightforward computation in \singular shows that $\dim (X,0)=1$ and that Sing$(X,0)=0$. Hence, $(X,0)$ is a one-dimensional \icis. \\

We have that $T\K_e f$ is the submodule generated by the elements $(2x,y), (2y,x), (2z,0)$ plus the submodule $(x^2+y^2+z^2, xy)\O_3^2$. With \singular\!, it can be verified that a $\C$-basis of $\m_3 \O_3^2 / T\K_e f$ is given by the vectors $(y,0), (0,y), (0,z)$. Thus, a miniversal deformation of $(X,0)$ is given by the map $F:(\C^3\times \C^3, 0)\rightarrow (\C^2\times \C^3,0)$ given by 
\begin{equation*}
    F(x,y,z,u,v,w)=(x^2+y^2+z^2+uy, xy+vy+wz, u,v,w). 
\end{equation*}
Furthermore, $\tau (X,0)=\codimKe (f)=5$.
\end{example}

\subsection{Milnor number of \textsc{icis} and Lê-Greuel's formula}
In this section, we review some key aspects of Milnor fibrations, which are necessary for defining the Milnor number of an \icis\!\!. We establish a remarkable result regarding the relationship between the Milnor and Tjurina numbers, which states that the former is always greater than or equal to the latter. Although this inequality may seem innocent, it turns out to connect two numerical invariants with completely different natures, resulting in an important outcome that motivates the Mond conjecture. Throughout this section, we focus on monogerms of \icis centered at the origin. \\

Let $X\subset \C^N$ be a representative of an analytic space at the origin that has isolated singularity. Let $g:(\C^N, 0)\rightarrow (\C^p,0)$ be a mapping such that $X=V(g_1, \ldots, g_p)$, and take a \textit{good representative} of $f$, which is a mapping of the form 
\begin{equation*}
    f:X\cap \overline{B_\epsilon}\cap f^{-1}(B_\eta) \rightarrow B_\eta,
\end{equation*}
where $0<\eta \ll \epsilon \ll 1$ (that is, take $\epsilon>0$ small enough, and once we have fixed $\epsilon$, then take $\eta>0$ small enough depending on $\epsilon$). In this case, if $s\in B_\eta$ is a regular value of $f$, we call $\overline{X_s}=f^{-1}(s)$ a \textit{Milnor fibre} of $X$. \\

The following central theorem was first noticed by Milnor for hypersurfaces, and then by Hamn for \icis\!:
\begin{theorem} If $(X,0)$ is an \icis of dimension $n$, then there exists a nonnegative integer $\mu$ such that, for any regular value $s\in B_\eta$, the Milnor fibre $\overline{X_s}$ has the homotopy type of a wedge of $\mu$ spheres $S^n$. 
\end{theorem}
\begin{definition} In the above conditions, we define the \textit{Milnor number} of $(X,0)$ as the number of spheres in the previous theorem, and we denote it as $\mu(X,0)$. 
\end{definition}
In what follows, we give some remarkable results concerning this definition. 
\begin{proposition} The Milnor number is invariant under isomorphism. That is, if $(X,0)$ and $(Y,0)$ are isomorphic \icis\!, then $\mu (X,0)=\mu(Y,0)$.   
\end{proposition}
A relevant algebraic formula is known for hypersurfaces. Namely, if $(X,0)=g^{-1}(0)$ for some function $g:(\C^N,0)\rightarrow (\C,0)$, then 
\begin{equation*}
    \mu (X,0) = \dim_{\C} \dfrac{\O_N}{J(f)}. 
\end{equation*}
For a generic \icis, determining the Milnor number is more intricate as there is no algebraic object that directly corresponds to its complex dimension. Nevertheless, there is a recursive method known as the Lê-Greuel formula that can be used to compute this invariant effectively. \\

Let $(X,0)$ be an \icis of dimension $n$, defined by as the fibre of a $\K$-finite mapping $g:(\C^{n+k},0)\rightarrow (\C^k,0)$. Since the discriminant $(D,0)$ of $f$ is known to be a hypersurface in $(\C^k,0)$, it follows that one can choose a line $L$ in such a way that $L\cap D={0}$. After a linear change of coordinates, we can assume that $L$ corresponds to the line $L$ given by $y_1=\ldots = y_{k-1}=0$. This implies that $g'=(g_1, \ldots, g_{k-1}):(\C^{n+k},0) \rightarrow (\C^{k-1},0)$ defines an \icis $(X',0)$ of dimension $n+1$ that contains $(X,0)$. In these terms, the following relationship between $(X,0)$ and $(X',0)$ holds:
\begin{theorem}[Lê-Greuel's formula] In the above situation, we have that 
\begin{equation*}
    \mu(X,0)+\mu (X',0)= \dim_\C \dfrac{\O_{n+k}}{I(g')+R(g)},
\end{equation*}
where $I(g')=(g_1, \ldots, g_{k-1})$ and $R(g)$ denotes the \textup{ramification ideal} of $g$, which is the ideal in $\O_{n+k}$ generated by the maximal minors of the jacobian matrix $dg$. 
\end{theorem}
This formula allows us to compute the Milnor number of any \icis by means of a recursive formula. Indeed, after a generic linear change of coordinates in $\C^k$ one can have that, for each $m\in \{1,\ldots, k\}$, the mapping $g^{(m)}:=(g_1, \ldots, g_m):(\C^{n+k},0)\rightarrow (\C^m,0)$ defines an \icis $(X^{(m)},0)$ of dimension $n+k-m$. It is straightforward to verify that 
\begin{theorem} In the above situation, we have that 
\begin{equation*}
    \mu (X,0) = \sum_{m=1}^k (-1)^{k-m} \dim_\C \dfrac{\O_{n+k}}{I(g^{(m-1)})+R(g^{(m)})}.
\end{equation*}
\end{theorem}
Notice that, if $k=1$, then we recover the well-known formula 
\begin{equation*}
    \mu (X,0)=\dim_\C \dfrac{\O_{n+1}}{J(g)}, 
\end{equation*}
where $R(g)$ coincides with the jacobian ideal $J(g)$. 
\begin{example} Let us compute the Milnor number of the \icis given by $(X,0)=V(x^2+y^2+z^2, xy)\subset (\C^3,0)$ of a previous example. In this case, $(X',0)=V(x^2+y^2+z^2)$, so $\mu (X',0)=1$. By Lê-Greuel's formula, we have that 
\begin{equation*}
    \mu (X,0) + \mu(X',0)=\dim_\C \dfrac{\O_3}{(x^2+y^2+z^2)+I_2\begin{pmatrix}
        2x & 2y & 2z \\ y & x & 0 
    \end{pmatrix}} =6,
\end{equation*}
where the previous dimension is computed with \singular\!. Hence, $\mu(X,0)=5$. 
\end{example}

Consider a hypersurface $(X,0)$ with an isolated singularity. We have defined two numerical invariants, namely the Milnor and Tjurina numbers. In this case, the invariants can be easily computed as
\begin{equation*}
    \mu(X,0)= \dim_\C\dfrac{\O_{n}}{J(g)}, \quad \tau (X,0)=\dim_\C \dfrac{\O_{n}}{(g)+ J(g)},
\end{equation*}
provided $g:(\C^n,0)\rightarrow (\C,0)$ is an equation of $(X,0)$. It is therefore obvious that $\mu(X,0)\geq \tau (X,0)$, since $J(g)\subset (g)+J(g)$. Furthermore, equality holds if and only if $g\in J(g)$, which can be shown to be equivalent to that $g$ is weighted homogeneous. Formally, if $g$ is weighted homogeneous, then $g\in J(g)$ trivially. The opposite direction, however, is nontrivial, and it was shown by Saito in 1971. \\

Although this observation is practically trivial, it has significant implications. It asserts that the number of spheres that appear after taking a smooth deformation $\overline{X_s}$ of $(X,0)$, which is $\mu(X,0)$, is always greater than or equal to the number of parameters required for a miniversal deformation of $(X,0)$, which is $\tau(X,0)$. Moreover, equality holds only if the germ $g$ is weighted homogeneous. \\

The spirit of this project is to extend the scope of this inequality to different frameworks. A natural question that arises in the context of this chapter is whether this inequality still hold if $(X,0)$ is an \icis\!\!. Although the result holds true in this case as well, its proof is rather nontrivial and required several more years to complete. The formal statement of the result is as follows:
\begin{theorem} Let $(X,0)$ be an \icis given by the vanishing of $g_1, \ldots, g_s$. Then, $\mu (X,0)\geq \tau (X,0)$ with equality if and only if each $g_i$ is weighted homogeneous. 
\end{theorem}
The fact that $\mu(X,0)\geq \tau (X,0)$ was shown by Looijenga and Steenbrick in 1985. In 1980, Greuel showed that weighted homogeneous \icis have $\mu = \tau$, and later in 2002, Vosegaard showed the opposite implication. \\

This kind of results are frequently referred as $\mu \geq \tau $ theorems. In the next chapter, we give the context to formalise a $\mu \geq \tau$ result that is still an open problem, which is the Mond conjecture. 

\newpage
\section{Image Milnor number and the Mond conjecture}
As we have commented on in the last chapter, if $(X,0)$ is a complete intersection with isolated singularity (that is, an \icis\!\!), we know that $\mu (X,0) \geq \tau (X,0)$. In other words, the number of spheres that appear in a smooth deformation of $(X,0)$ given by a Milnor fibre is greater than the number of parameters of a miniversal deformation of $(X,0)$. Moreover, equality holds if and only if the germ $(X,0)$ is weighted homogeneous, meaning that all its equations are weighted homogeneous. In this chapter, we focus on a different type of analytic germs, namely hypersurfaces that can be described as the image of an $\A$-finite mapping $f:(\mathbb{C}^n,S)\rightarrow (\mathbb{C}^{n+1},0)$. These hypersurfaces do not have an isolated singularity, and therefore their analysis must be approached differently. We investigate the natural context in which a $\mu\geq \tau$ result can be formulated, which is known as the Mond conjecture.

\subsection{The image Milnor number}
In this section, we present one of the main concepts of this project, namely the image Milnor number. This concept plays the role of the Milnor number of a hypersurface, as it measures the number of shperes that arise after a perturbation of the set.  Firstly, we present some concepts that are related with the image Milnor number, and direct our focus towards a concept introduced by Mather: the nice dimensions. 
\begin{definition}
    We say that $(n,p)$ are \textit{nice dimensions} if the set of all stable mappings $f:U\subset \C^n \rightarrow V\subset \C^p$ is dense in $C^\infty_\text{pr}(U,V)$.
\end{definition}
The particular topology that is given to measure the density, as well as the given space of mappings, is not relevant for our purposes. Hence, we omit the details and refer them to \cite{juanjo}, sec. 5.2.2. The intuition that is relevant in this definition is that, in the nice dimensions, any given mapping can be approximated as nicely as 
\begin{multicols}{2}
neccesary by a stable mapping. 
The region of the nice dimensions is given in figure 1, where $(n,p)$ are in the nice dimensions if the corresponding point is contained in the region that lies in the left or below the polygon. It is of particular interest for our purposes to notice that $(n,n+1)$ are nice dimensions if and only if $1\leq n \leq 14$. 

\begin{center}
        \centering
        \includegraphics[width=0.4\textwidth]{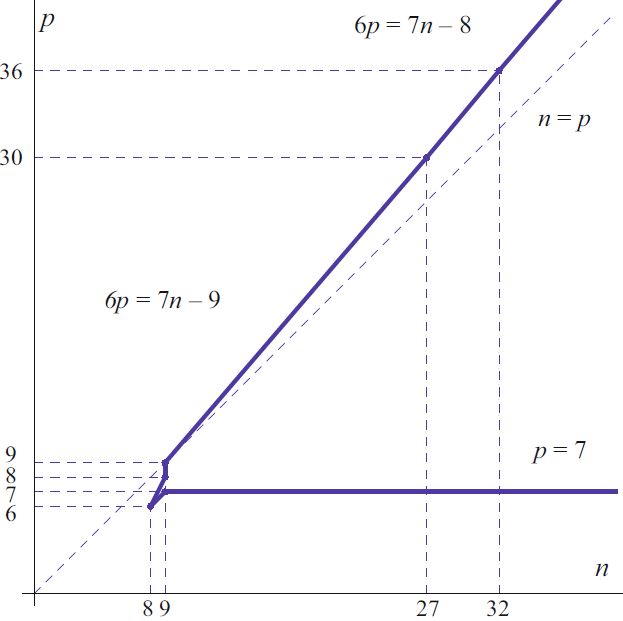}
        \textbf{Figure 1:} Mather's nice dimensions
\end{center}
\end{multicols}
As a hypersurface that lacks isolated singularity does not possess a clearly defined Milnor number, we must identify an equivalent numerical invariant that is applicable within this new context. To this end, we introduce the analogous concept of the Milnor fibration in this setting.
\begin{definition} A \textit{stabilisation} of a mapping $f:(\C^n,S)\rightarrow (\C^p,0)$ is a $1$-parameter unfolding $F(x,t)=(f_t(x),t)$, for each $t\neq 0$, the mapping $f_t$ is locally stable. 
\end{definition}
Thus, a stabilisation of a mapping $f$ is a deformation $f_t$ of the given mapping in such a way that $f_0=f$ and that, for each $t\neq 0$, the deformation $f_t$ is stable. Hence, $f_t$ has the simplest possible singularities for $t\neq 0$. The following theorem assures the existence of stabilisations under certain conditions.
\begin{proposition} If an $\A$-finite mapping $f:(\C^n,S)\rightarrow (\C^p,0)$ satisfies that $(n,p)$ are nice dimensions, then $f$ admits a stabilisation. 
\end{proposition}
For a reference, see proposition 5.5 of \cite{juanjo}. 
\begin{remark}
    It is worth noting that there exist other conditions under which the existence of a stabilisation can be guaranteed. For instance, any map-germ of corank one admits a stabilisation. However, in this project, we focus on ensuring the existence of a stabilisation when the dimensions $(n,p)$ are nice dimensions.
\end{remark}

Let us direct our focus to hypersurfaces defined as the image of a mapping. Consider an $\A$-finite mapping $f:(\C^n,S)\rightarrow (\C^{n+1},0)$, where our interest lies in studying the image $(X,0)$ of $f$. It can be verified that any $\A$-finite mapping is finite and generically one-to-one, thereby the set-germ of image $(X,0)$ is well-defined and has dimension $n$. Thus, $(X,0)$ is a hypersurface in $(\C^{n+1},0)$. In the field of analytic geometry, it is customary to refer to the restriction $f:(\C^n,0)\rightarrow (X,0)$ as the \textit{normalisation} of $(X,0)$, which essentially implies that $f$ acts as a parametrisation for $(X,0)$. To delve deeper into the geometry of $(X,0)$, we leverage the machinery of singularities of mappings. \\

As we have commented on before, provided $(n,n+1)$ are nice dimensions, any $\A$-finite mapping $f:(\C^n,S)\rightarrow (\C^{n+1},0)$ admits a stabilisation $f_t$. Let us denote by $X_t$ to the image of $f_t$ for $t\neq 0$. As we have already mentioned, $X_t$ has the simplest possible singularities, owing to the stability of $f_t$. We refer to $X_t\cap \overline{B_\epsilon}$ as a \textit{disentanglement} of $f$, where $\overline{B_\epsilon}$ denotes the closed ball of center $0$ and small enough radius $\epsilon$. \\

The following theorem establishes that this disentanglement $X_t$ acts as a Milnor fibre if we were in the case of a hypersurface having an isolated singularity. More specifically,
\begin{theorem} Assume that $f:(\C^n,S)\rightarrow (\C^{n+1},0)$ is an $\A$-finite mapping that admits a stabilisation $f_t$. Then, there exists a non-negative integer $\mu$ such that any disentanglement of $f$ has the homotopy type of a wedge of $\mu$ spheres of dimension $n$. 
\end{theorem}
In the present setting, we introduce the notion of the \textit{image Milnor number} of $f$, which counts the number of spheres that appear in the topology of a disentanglement of $f$, denoted by $\mu_I(f)$. This concept is analogous to the Milnor number of a hypersurface with an isolated singularity or an \icis. In other words, it quantifies the number of spheres that emerge following a ``nice'' perturbation of the image set of $f$. \\

Altough it may seem, by the way it is denoted, that the image Milnor number $\mu_I(f)$ or the codimension $\codimAe(f)$ depends on the mapping $f$ that parametrises its image $(X,0)$, it turns out that both are analytical invariants of $(X,0)$. 
\begin{proposition} If $f_1,f_2:(\C^n,S)\rightarrow (\C^{n+1},0)$ are $\A$-finite mappings whose images $(X_1, 0), (X_2,0)$ are induced isomorphic, then $f_1$ and $f_2$ are $\A$-equivalent.
\end{proposition}
\begin{proof} If $\psi:(\C^{n+1},0)\rightarrow (\C^{n+1},0)$ is the isomorphism that takes $X_1$ to $X_2$, then $\psi\circ f_1$ and $f_2$ are both finite and generically one-to-one mappings that parametrise the same image $(X_2,0)$. As the normalisation map is unique up to a change of coordinates in the source and both $\psi \circ f_1$ and $f_2$ are normalisations of $X_2$, it follows that there exists an isomorphism $\phi:(\C^n,S)\rightarrow (\C^n,S)$ such that $\psi \circ f_1 = f_2 \circ \phi$. In particular, $f_1$ and $f_2$ are $\A$-equivalent. 
\end{proof}
Now, it is clear that $\A$-equivalent mappings have the same image Milnor number and codimension. Hence, $\mu_I(f)$ and $\codimAe(f)$ are, in fact, invariants of the analytical type of the image $(X,0)$ of $f$. This makes clear that $\mu_I(f)$ and $\codimAe(f)$ are both invariants of $(X,0)$, and will be shown to extend the well-known definitions of Milnor and Tjurina numbers in the case of isolated singularity. 
\subsection{The Mond conjecture}
In this section, we present the formal statement of the Mond conjecture, which is a problem of the form $\mu\geq \tau$. This means that the problem is analogous to the well-known inequality $\mu(X,0)\geq \tau (X,0)$ for hypersurfaces with isolated singularity. In that context, $\mu(X,0)$ is the number of spheres that appear in the topology of the Milnor fiber of $(X,0)$, which is a smooth approximation of $(X,0)$, and $\tau(X,0)$ is the number of parameters of a miniversal deformation of $(X,0)$. Moreover, the equality $\mu=\tau$ holds if and only if $(X,0)$ is weighted homogeneous.

In our current setting, we analyse $\A$-finite mappings $f:(\C^n,S)\rightarrow (\C^{n+1},0)$, for which we have established that $\mu_I(f)$ measures the number of spheres that appear in the topology of the disentanglement of $f$. Furthermore, $\codimAe(f)$ is known to represent the number of parameters that an $\A$-miniversal unfolding of $f$ must have. Therefore, the question becomes evident: is it true that
$$\mu_I(f)\geq \codimAe (f)?$$
\begin{conjecture}[Mond conjecture] Let $f:(\C^n,S)\rightarrow (\C^{n+1},0)$ be an $\A$-finite mapping, where $(n,n+1)$ are nice dimensions. Then, $\mu_I(f)\geq \codimAe (f)$ with equality if $f$ is weighted homogeneous.
\end{conjecture}
This result therefore intends to extend the scope of the well-known analogous result for hypersurfaces with isolated singularity. Indeed, this is the analogous statement for the hypersurfaces $(X,0)$ that can be described as the image of an $\A$-finite mapping $f:(\C^n,S)\rightarrow (\C^{n+1},0)$, and thus $(X,0)$ has parametric equations given by $f$. \\

At present, the only values of $n$ for which the conjecture is known to be true is for $n=1$ and $n=2$. In other words, the inequality holds for analytic spaces $(X,0)$ that are parametric curve singularities in $(\C^2,0)$ and parametric surface singularities in $(\C^3,0)$. \\

Nevertheless, it is not generally true that any analytic space $(X,0)$ of dimension $n$ can be described as the image of an $\A$-finite mapping $f:(\C^n,S)\rightarrow (\C^{n+1},0)$. In analytic geometry, these germs are referred to as the ones with ``smooth normalisation'', meaning that they admit a normalisation $f:(\C^n,S)\rightarrow (X,0)$ to the smooth space $(\C^n,S)$. Our contribution will be to expand this result to include surfaces that do not necessarily have a smooth normalization. Specifically, we extend the Mond conjecture to surfaces whose normalization is an \textsc{icis}, and where the normalisation mapping is $\A$-finite. In order to do so, we need to study and develop techniques for mappings defined on \textsc{icis}. This will be the focus of the following chapter.

\newpage
\thispagestyle{empty}
\textcolor{white}{blank page}
\newpage
\section{Singularities of mappings defined on \icis}
In the first three chapters, we examined singularities of mappings of the form $f:(\C^n,S)\rightarrow (\C^p,0)$, as well as the properties of a particular case of analytic spaces in $(\C^n,S)$, which are known as \icis. In fact, these are the two main objects of study in singularity theory: mappings and sets. We have shown that the study of deformations in both of them is closely related. In this second part of this project, we explore mappings whose domains are allowed to have singularities. Specifically, we investigate mappings $f:(X,0)\rightarrow (\C^p,0)$ defined on \icis $(X,0)\subset (\C^N,0)$. In their paper \cite{Mond-Montaldi}, Mond and Montaldi established a basic framework for studying such mappings, where deformations are allowed to deform both the mapping $f$ and the domain $(X,0)$. We rely mainly on this article by Mond and Montaldi, as well as the clear exposition provided by Giménez Conejero and Nuño-Ballesteros in \cite{juanjo-roberto}. The main objective of this chapter is to introduce an analogous definition of the image Milnor number in this context, as well as to state a generalised version of the Mond conjecture for these mappings. \\

Throughout this chapter, we use $(X,S)\subset (\C^N,S)$ to denote a multi-germ of an \icis, and $f:(X,S)\rightarrow (\C^p,0)$ denotes a holomorphic mapping, in the sense that it can be extended to a holomorphic mapping $\tilde{f}:(\C^N,S)\rightarrow (\C^p,0)$. We commonly denote these mappings as $(X,f)$ to highlight their dependence on the domain.
\subsection{Deformations of mappings defined on \icis}
In this section, we present the definitions and results regarding deformation theory of these mappings, and provide a definition for the codimension for mappings defined on \icis\!\!. We follow the same structure as in section 1.1, where $\A$-equivalence of mappings was studied. 
\begin{definition} Let $f,g:(X,S)\rightarrow (\C^p,0)$ be germs of mappings. We say that $f$ and $g$ are $\A$\textit{-equivalent} if there exist isomorphisms $\phi:(X,S)\rightarrow (X,S)$ and $\psi:(\C^p,0)\rightarrow (\C^p,0)$ such that the diagram
    \begin{center}
        \begin{tikzcd}
          (X,S) \arrow[r, "f"] \arrow[d, "\phi"]
          & (\C^p,0) \arrow[d, "\psi"] \\
          (X,S) \arrow[r, "g"]
          & (\C^p,0)
          \end{tikzcd}
      \end{center}
    commutes. 
\end{definition}
Hence, singularities of mappings with \icis in the source will be studied up to $\A$-equivalence, in the same way as it was done in case that $(X,S)$ is smooth. 
\begin{definition} An $r$\textit{-parameter unfolding} of $f:(X,S)\rightarrow (\C^p,0)$ is a mapping $F:(\X,S')\rightarrow (\C^p\times \C^r,0)$ together with a projection $\pi :(\X,S')\rightarrow (\C^r,0)$ such that $(X,S)$ is isomorphic to the fibre $(\pi^{-1}(0),S')$, and if $i:(X,S)\rightarrow (\X,S')$ denotes the natural embedding associated with the given isomorphism, and $j:(\C^p,0)\rightarrow (\C^p\times \C^r,0)$ is the mapping $j(y)=(y,0)$, then the diagram 
\begin{center}
        \begin{tikzcd}
	{(X,S)} & {(\C^p,0)} \\
	{(\X,S')} & {(\C^p\times \C^r,0)} \\
	& {(\C^r,0)}
	\arrow["f", from=1-1, to=1-2]
	\arrow["i"', hook, from=1-1, to=2-1]
	\arrow["j", hook, from=1-2, to=2-2]
	\arrow["F", from=2-1, to=2-2]
	\arrow["\pi"', two heads, from=2-1, to=3-2]
	\arrow["{\pi_2}", two heads, from=2-2, to=3-2]
\end{tikzcd}
      \end{center}
commutes. We ofter shorten the notation to say that $(\X, \pi , F)$ is the unfolding.  
\end{definition}
\begin{remark} In order to verify that the given definition of unfloding matches with the intuition that lies behind the concept, let us check that this definition indeed recovers the given for mapping defined on a smooth source $(X,S)=(\C^n,S)$. Indeed, in this case, $(\X,S')=(\C^n\times \C^r,S\times 0)$ has a natural projection $\pi:(\C^n\times \C^r,S\times 0)\rightarrow (\C^r,0)$ for which $(\C^n,S)\cong \pi^{-1}(0)=(\C^n\times 0, S\times 0)$. The commutativity of the square in the diagram implies that $F(x,0)=(f(x),0)$ for each $x$, and the commutativity of the triangle forces $F$ to have the form $F(x,u)=(f_u(x), u)$. Hence, $F$ is an unfolding in the standard case. 
\end{remark}
\begin{remark} As in the definition of a deformation of an \icis, notice that the concept depends on $(X,S)$ up to isomorphism. Thus, after applying an isomorphism in the source, one can assume that $(X,S)$ \textit{coincides} with the fibre $(\pi^{-1}(0),S')$, where $S=S'$. Hence, the given embedding becomes an inclusion. With this identification, the intuition of the concept becomes clearer. The domain of an unfolding is a germ $(\X,S)$ that contains $(X,S)$ as the fibre $\pi^{-1}(0)$, and that contains its deformations $X_u=\pi^{-1}(t)$. Furthermore, the mappings $f_t=\pi_1\circ F |_{X_t}:X_t\rightarrow \C^p$ serve as \textit{perturbations} of $f_0=f$, where both the domain and the mapping itself are being deformed.
\end{remark}
\begin{example}A particular kind of unfolding that arises is where $\X = X\times \C^r$. In this case, the unfolding only permits deformations on the mapping $f$, but not on the source $(X,S)$. 
\end{example} 
\begin{definition} Two $r$-parameter unfoldings $(\X, \pi , F)$ and $(\X',\pi',F')$ of the same mapping $f:(X,S)\rightarrow (\C^p,0)$ are said to be $\A$\textit{-equivalent} if there exist germs of isomorphisms $\Phi:(\X,S)\rightarrow (\X',S')$ and $\Psi:(\C^p\times\C^r, 0)\rightarrow (\C^p\times\C^r, 0)$ such that the diagram 
\begin{equation*}
    \begin{tikzcd}
	{(\X,S)} & {(\C^p\times \C^r,0)} \\
	{(\X',S')} & {(\C^p\times \C^r,0)}
	\arrow["F", from=1-1, to=1-2]
	\arrow["\Phi"', from=1-1, to=2-1]
	\arrow["\Psi", from=1-2, to=2-2]
	\arrow["{F'}", from=2-1, to=2-2]
\end{tikzcd}
\end{equation*}
commutes, and where
\begin{enumerate}
    \item $\Psi$ is an unfolding of the identity $\id_{(\C^p,0)}$,
    \item $\Phi|_{(X,S)}=\id_{(X,S)}$, and
    \item $\pi'\circ \Phi = \pi$, \textit{i.e.}, $\Phi$ sends the fibers $X_t$ in $\X$ to the fibers $X'_t$ in $\X'$ for each $t$. 
\end{enumerate}
\begin{example}\label{example:unfolding}  Let us consider the germ $(X,0)\subset (\C^2,0)$ given by the equation $x^2-y^3=0$, and the mapping $f:(X,0)\rightarrow (\C,0)$ given by $f(x,y)=y^2$. A possible unfolding of $f$ is the one given over the space $(\X,0)\subset (\C^3,0)$ described by the equation $x^2-y^3-ty(2y+t)=0$, and $F:(\X,0)\rightarrow (\C^2,0)$ is the mapping given by $F(x,y,t)=(y^2+ty,t)$. Then, if $\pi: (\X,0)\rightarrow (\C,0)$ is given by $\pi(x,y,t)=t$, it is clear that $\pi^{-1}(0)=(X,0)$ makes $F$ an unfolding of $f$. 
    \begin{figure}[h]
        \centering
        \includegraphics[width=0.6\textwidth]{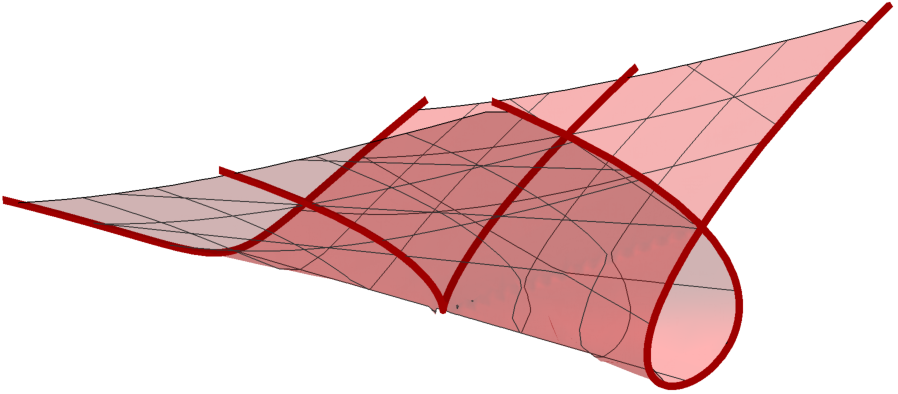}
        \caption{Real picture of the space $\X$ in the example \ref{example:unfolding}. The red fibre in the middle corresponds to $X=\pi^{-1}(0)$ and the other ones correspond to $\pi^{-1}(-1)$ and $\pi^{-1}(1)$. }
        \label{fig:my_label}
    \end{figure}
\end{example}
\end{definition}
The simplest kind of unfoldings are the ones that do not deform the mapping in any sense:
\begin{definition} An unfolding of $f:(X,S)\rightarrow (\C^p,0)$ is $\A$\textit{-trivial} if is is $\A$-equivalent to the constant unfolding $F:(X\times \C^r,S\times 0)\rightarrow (\C^p\times \C^r,0)$ given by $F(x,u)=(f(x),u)$. 
\end{definition}
In these terms, we say that $(X,f)$ is \textit{stable} provided every unfolding is $\A$-trivial. In other case, we say that the mapping $(X,f)$ is \textit{unstable}. 

\begin{definition} If $(\X, \pi , F)$ is an $r$-parameter unfolding of $(X,f)$ and $\rho:(\C^s, 0)\rightarrow (\C^r,0)$ is a holomorphic mapping, then $\rho$ induces an $s$-parameter unfolding $(\X', \pi', F')$ by a \textit{base change} through the fibered product of $F$ and $\id_{(\C^p,0)}\times \rho$ through the commutative diagram
\begin{equation*}
    \begin{tikzcd}
	{\X'=\X \times_{\C^p\times\C^r} (\C^p\times \C^s)} & {(\C^p\times \C^s,0)} \\
	{(\X',S')} & {(\C^p\times \C^r,0),}
	\arrow["{F'}", from=1-1, to=1-2]
	\arrow[from=1-1, to=2-1]
	\arrow["{\id_{\C^p}\times \rho }", from=1-2, to=2-2]
	\arrow["F", from=2-1, to=2-2]
\end{tikzcd}
\end{equation*}
where the germ notation is ommited for simplicity. In this case, we say that $(\X',\pi', F')$ is a pull-back of the unfolding $(\X, \pi , F)$ via $\rho$. 
\end{definition}
\begin{definition} An unfolding $(\X, \pi , F)$ is $\A$-versal provided that any other unfolding is $\A$-equivalent to a pull-back of $(\X, \pi , F)$. We say that $(\X, \pi , F)$ is $\A$\textit{-miniversal} whenever it is $\A$-versal with the property that no other $\A$-versal unfoldings exist with less parameters. 
\end{definition}
The intuition of this concept is, as in the smooth case, that an $\A$-versal unfolding does reflect all possible deformations of $(X,f)$. moreover, an $\A$-miniversal unfolding encodes all these deformations without redundancies.  
\begin{definition} The number of parameters of an $\A$-miniversal unfolding of $(X,f)$ is denoted as its $\A_e$\textit{-codimension} if it is finite, and it is denoted as $\codimAe (X,f)$. If this number does not exist, we set the codimension to be infinite. If $(X,f)$ has finite codimension, we say that $(X,f)$ is $\A$-finite. 
\end{definition}
It turns out that Mather-Gaffney's criterion can be adapted to this context to yield
\begin{proposition} A mapping $(X,f)$ is $\A$-finite if and only if it has isolated instability. 
\end{proposition}

A relevant remark that has to be done here is that $\A$-finite mappings are finite and generically one-to-one, as in the case of mappings with smooth source. This will be an important fact to take into account in the following chapters. \\

It is worth to mention that there exists an analogous version of the codimension of a mapping with non-smooth domain, where only deformations that do not deform the source are allowed. This codimension is often denoted as $\codimAe(f)$, where the domain does not appear as in $\codimAe(X,f)$. In fact, it is not a surprising fact that both invariants are related through the equality 
\begin{equation*}
    \codimAe(X,f)=\codimAe(f)+\tau(X). 
\end{equation*}
As it will be shown in later examples, the correct version of the generalised conjecture will involve $\codimAe(X,f)$. However, we suspect that another version of the conjecture could be stated in terms of $\codimAe(f)$ by using another definition of $\mu_I$. 

\subsection{The image Milnor number and the extended Mond conjecture}
In this section, we follow a similar procedure as we did in Chapter 3 to provide a rigorous definition for the image Milnor number and to settle the Mond conjecture in this general framework. To do so, we first need to define the concept of stabilization in this context. \\

In what follows, we let $(X,S)\subset (\C^N,S)$ be an \icis of dimension $n$ and $f:(X,S)\rightarrow (\C^p,0)$ an an $\A$-finite mapping. 
\begin{definition} A \textit{stabilisation} of $f$ is a one-parameter unfolding $(\X,\pi, F)$ where its associated perturbations $f_t:X_t\rightarrow \C^p$ are locally stable for $t\neq 0$. 
\end{definition}
Similar to the case of mappings with smooth domains, the existence of stabilizations can be ensured in Mather's nice dimensions.
\begin{proposition} If $(n,p)$ are nice dimensions, then $f$ admits a stabilisation. 
\end{proposition}
Hence, under the above conditions, $f$ admits a stabilisation $f_t:X_t\rightarrow \C^{n+1}$. In this case, we say that, for $t\neq 0$, the image $Y_t\cap \overline{B_\epsilon}=f_t(X_t)\cap \overline{B_\epsilon}$ is a \textit{disentanglement} of $(X,f)$, provided $\epsilon>0$ is small enough.  
\begin{theorem} Let $(X,S)$ be an $n$-dimensional \icis and let $f:(X,S)\rightarrow (\C^{n+1},0)$ be an $\A$-finite mapping that admits a stabilisation. Assume that $(n,n+1)$ are nice dimensions. Then, there exsits a non-negative integer $\mu$ such that any disentanglement of $(X,f)$ has the homotopy type of a wedge of $\mu$ spheres of dimension $n$. 
\end{theorem}

The proof of both results is easily deduced from the original proof in case that $f$ is defined over a smooth space, with minor adaptations. \\

In this context, we define the \textit{image Milnor number} of $(X,f)$ as the number of spheres that appear in the topology of a disentaglement, and we denote it as $\mu_I(X,f)$. 

\begin{figure}[h]
    \centering
    \includegraphics[width=0.6\textwidth]{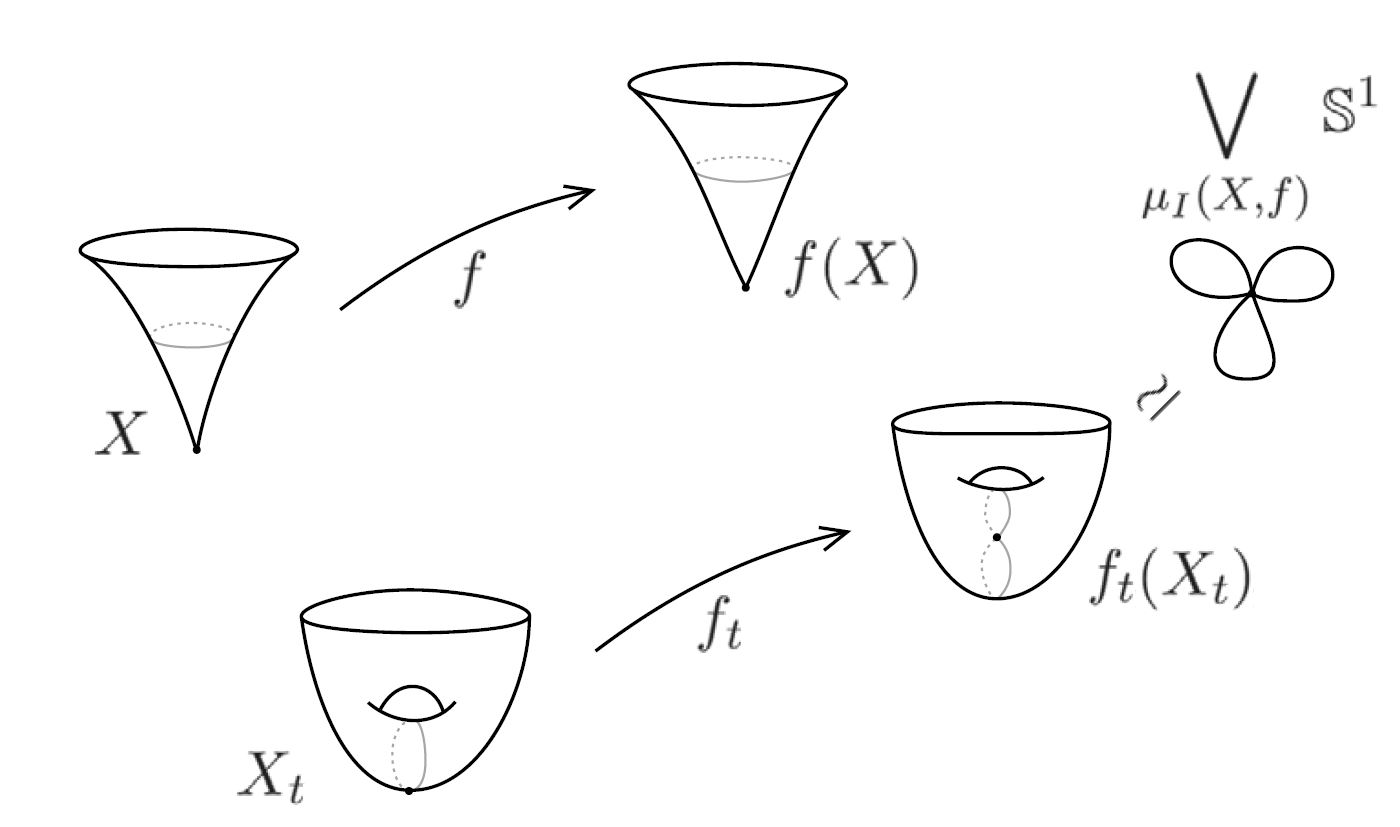}
    \caption{Representation of the disentanglement of a mapping defined on an ICIS. Image extracted from the
article \cite{roberto} by R. Giménez Conejero and J.J. Nuño-Ballesteros.}
    \label{fig:enter-label}
\end{figure}

Now that we have defined the image Milnor number and the codimension of a mapping $(X,f)$, we have all the necessary components to state an analogous version of the Mond conjecture.
\begin{conjecture} [Extended Mond conjecture] Let $(X,S)$ be an $n$-dimensional \icis and $f:(X,S)\rightarrow (\C^{n+1},0)$ be an $\A$-finite mapping, where $(n,n+1)$ are nice dimensions. Then $\mu_I(X,f)\geq \codimAe (X,f)$, with equality if $(X,f)$ is weighted homogeneous. 
\end{conjecture}
Recall that the standard conjecture is currently known to be true for $n=1,2$. In this general setting, however, it has only been verified for $n=1$ by Nuño and Henrique in \cite{juanjo-henrique}. One of the main objectives of this project is to show that the conjecture holds for $n=2$ in this extended setting. In order to do so, we define in the following chapter a module that is able to measure the image Milnor number. \\

As we have mentioned before, $\mu_I(X,f)$ and $\codimAe(X,f)$ are both analogous invariants of the Milnor and Tjurina numbers for a particular kind of hypersurfaces that do not have isolated singularity. In this sense, the hypersurfaces that we deal with are the ones that have normalisation map with isolated instability (that is, $\A$-finite mappings) where the associated normal space is an \textsc{icis}. Therefore, this project extends the scope in which the $\mu\geq \tau$ can be defined and proved. Some examples of the conjecture will be shown in the subsequent chapters. 

\newpage
\section{A Jacobian module for disentanglements with smooth source}
The main impediment in proving the Mond conjecture lies in the fact that the image Milnor number lacks an algebraic expression for computation. To overcome this challenge, in this chapter we construct a module $M(g)$ that has the image Milnor number as its dimension under certain conditions. Consequently, we show that the Mond conjecture holds provided this condition is satisfied. To achieve this, we follow the methods outlined in \cite{BobadillaNunoPenafort} with the intention of adapting them for mappings defined on \icis in the following chapter.
\subsection{The module \textit{M(g)} for mappings defined on smooth sources}
In this initial section, the fundamental results needed for the proper construction of the module $M(g)$ are presented, with the aim of measuring the image Milnor number. We follow the construction provided in \cite{BobadillaNunoPenafort}. \\

Let $f:(\C^n,S)\rightarrow (\C^{n+1},0)$ be an $\A$-finite germ of mapping. Let $(Y,0)$ be the image of $f$ in $(\C^{n+1},0)$. The restriction $\overline{f}:(\C^n,S)\rightarrow (Y,0)$ is therefore the normalization map, and hence $\overline{f}^*:\O_{Y,0}\rightarrow \O_n$ is a monomorphism. Hence, $\O_{Y,0}$ can be seen as a subring of $\O_n$. In this case, the diagram
\begin{equation*}
\begin{tikzcd}
\O_{n+1} \arrow[r, "f^*"] \arrow[dr,"\pi",two heads]
&  \O_n \\
& \O_{Y,0}\arrow["i",u,hook]
\end{tikzcd}
\end{equation*}
commutes, where $\pi$ is the epimorphism induced by the inclusion of $(Y,0)$ in $(\C^{n+1},0)$. We consider both $\O_{Y,0}$ and $\O_n$ as $\O_{n+1}$-modules via the corresponding morphisms. Since every $\A$-finite mapping is finite, it follows that $(Y,0)$ is a hypersurface. Take $g:(\C^{n+1},0)\rightarrow (\C,0)$ a reduced equation of $(Y,0)$, and denote $J(g)\subset \O_{n+1}$ to the jacobian ideal of $g$. \\

Let us denote by $C(f)$ to the conductor ideal of $\O_{Y,0}$ in $\O_n$, which is given by 
\begin{equation*}
    C(f)=\{h\in \O_{Y,0} : h \cdot \O_n \subset \O_{Y,0} \}.
\end{equation*}
The conductor can be shown to be the largest ideal of $\O_{Y,0}$ which is also an ideal in $\O_n$. The following result from Piene \cite{Piene} gives a method to easily determine the conductor ideal. 

\begin{lemma}\label{lemma:piene}
    There exists a unique element $\lambda \in \O_n$ such that 
    \begin{equation*}
        \dfrac{\partial g}{\partial y_k}\circ f=(-1)^k\cdot \lambda \cdot\det (df_1, \ldots, df_{k-1}, df_{k+1}, \ldots, df_{n+1}),
    \end{equation*}
    for every $k\in \{1,\ldots, n+1\}$. Moreover, the ideal $C(f)$ is generated by $\lambda$. 
\end{lemma}

In particular, this lemma implies that $J(g)\cdot \O_n\subset C(f)$. Furthermore, $J(g)\subset (f^*)^{-1}(C(f))=\pi^{-1}(C(f))$. 
\begin{remark}\label{remark:F1_determinantal}
    It follows from the proof of Theorem 3.4 of \cite{Mond-Pellikaan} that $\pi^{-1}(C(f))$ coincides with the first Fitting ideal of $\O_n$ as an $\O_{n+1}$-module via $f^*$, which will be denoted as $\F_1(f)$. In fact, this theorem states that $\O_{n+1}/\F_1(f)$ is a determinantal ring of dimension $n-1$. Furthermore, it follows by proposition 1.5 of \cite{Mond-Pellikaan} that $V(\F_1(f))$ consists of the set-germ of points $y\in \C^{n+1}$ such that, either $y=f(x)$ for some non-immersive point $x\in \C^n$, or $y=f(x)=f(x')$, for some $x\neq x'$. Hence, $V(\F_1(f))$ is the singular locus of $(Y,0)$. Moreover, these results hold in the general setting of mappings defined on \icis, which will play a key role later.
\end{remark}
\begin{definition} The map $f^*$ naturally induces an epimorphism of $\O_{n+1}$-modules 
\begin{equation*}
    \dfrac{\F_1(f)}{J(g)}\rightarrow \dfrac{C(f)}{J(g)\cdot \O_{X,0}}.
\end{equation*}
We define $M(g)$ as the $\O_{n+1}$-module given as the kernel of the previous morphism. 
\end{definition}
\begin{remark} Naturally, the module $M(g)$ is determined by the exact sequence 
\begin{equation*}
    0\rightarrow M(g) \rightarrow \dfrac{\F_1(f)}{J(g)}\rightarrow \dfrac{C(f)}{J(g)\cdot \O_{X,0}} \rightarrow 0,
\end{equation*}
where the morphisms are the obvious ones. 
\end{remark}
A straightforward consequence of the definition is this useful formula to compute the module $M(g)$. 
\begin{proposition}\label{prop:formulaMg}
\begin{equation*}
    M(g)=\dfrac{(f^*)^{-1}(J(g)\cdot \O_{X,0})}{J(g)}.
\end{equation*}
\end{proposition}
\begin{proof} It is clear by definition that 
\begin{equation*}
    M(g)=\dfrac{(f^*)^{-1}(J(g)\cdot \O_{X,0})\cap \F_1(f)}{J(g)},
\end{equation*}
but $J(g)\cdot \O_{X,0}\subset C(f)$, which forces $(f^*)^{-1}(J(g)\cdot \O_{X,0})\subset \F_1(f)$. The claim follows. 
\end{proof}

The following natural step is to check whether the dimension of $M(g)$ upper bounds the codimension of the mapping $f$. In the original article \cite{BobadillaNunoPenafort}, the given proof does not generalise properly to the case that will be presented in the subsequent sections. Nevertheless, as we need this result to show the general version for mappings defined on \icis, we sketch its proof. One of the methods that can be followed to prove this is by applying the following formula for the codimension:
\begin{lemma}[Proposition 2.1 of \cite{vanishingcycles}] If $f:(\C^n,S)\rightarrow (\C^{n+1},0)$ is an $\A$-finite mapping and $n\geq 2$, then 
\begin{equation*}
    \codimAe (f) = \dim_\C \dfrac{J(g)\cdot \O_n}{J(g) \cdot \O_{Y,0}}. 
\end{equation*}
\end{lemma}
\begin{proposition} If $n\geq 2$ and $K(g)=((g)+J(g))/J(g)$, then 
$$\dim_\C M(g)=\dim_\C K(g) + \codimAe (f).$$
\end{proposition}
\begin{proof}
    In virtue of the previous lemma, it is enough to check that the sequence
    \begin{equation*}
        0 \rightarrow K(g) \rightarrow M(g) \rightarrow \dfrac{J(g)\cdot \O_n}{J(g)\cdot \O_{Y,0}}\rightarrow 0
    \end{equation*}
    is short exact, and then apply lengths. Consider the commutative diagram
    \begin{equation*}
        \begin{tikzcd}
	& 0 & 0 & 0 \\
	& {K(g)} & {M(g)} & {\frac{J(g)\cdot \O_n}{J(g)\cdot \O_{Y,0}}} \\
	0 & {K(g)} & {\frac{C(f)}{J(g)}} & {\frac{C(f)}{J(g)\cdot \O_{Y,0}}} & 0 \\
	& 0 & {\frac{C(f)}{J(g)\cdot \O_n}} & {\frac{C(f)}{J(g)\cdot \O_n}} & 0 \\
	&& 0 & 0
	\arrow[from=1-4, to=2-4]
	\arrow[from=1-3, to=2-3]
	\arrow[from=1-2, to=2-2]
	\arrow[from=2-2, to=2-3]
	\arrow[from=2-3, to=2-4]
	\arrow[from=3-1, to=3-2]
	\arrow[from=3-2, to=3-3]
	\arrow[from=3-3, to=3-4]
	\arrow[from=3-4, to=3-5]
	\arrow["{\mu_3}"', from=2-4, to=3-4]
	\arrow["{\mu_2}"', from=2-3, to=3-3]
	\arrow["{\mu_1}"', from=2-2, to=3-2]
	\arrow["{\lambda_1}"', from=3-2, to=4-2]
	\arrow["{\lambda_2}"', from=3-3, to=4-3]
	\arrow["{\lambda_3}"', from=3-4, to=4-4]
	\arrow[from=4-4, to=5-4]
	\arrow[from=4-3, to=5-3]
	\arrow[from=4-3, to=4-4]
	\arrow[from=4-2, to=4-3]
	\arrow[from=4-4, to=4-5]
\end{tikzcd}
    \end{equation*}
    where $\lambda_i$ and $\mu_i$ are the natural mappings. Since the columns and the second and third rows are exact, it follows by Snake Lemma that the first row is exact. 
\end{proof}
\begin{remark} Notice that, if $f$ is stable, then $\dim_\C M(g)=0$. Indeed, in this case $\codimAe (f)=0$. As it is shown in section 7.4 of \cite{juanjo}, all stable mappings have a quasi-homogeneous equation $g$, and hence $K(g)=0$. Hence, $M(g)=0$ if and only if $f$ is stable and $g\in J(g)$. 
\end{remark}

\subsection{The relative version for unfoldings of mappings with smooth source}\label{section:5_2}
In this subsection, we examine the behavior of the module $M(g)$ under deformations. To achieve this, we introduce a relative version of this module for unfoldings, which reduces to the original module when the parameters are zero. To ensure clarity, we have placed some technical details concerning tensor products and flatness in Appendix \ref{appendix:1}. We encourage the reader to review these concepts prior to reading this section. \\

Let $F$ be an $r$-parameter unfolding of $f:(\C^n,S)\rightarrow (\C^{n+1},0)$. Then, the diagram 
\begin{equation*}
        \begin{tikzcd}
          (\C^{n+r},S) \arrow[r, "F"] 
          & (\C^r\times \C^{n+1},0)  \\
          (\C^n,S) \arrow[u, "i"] \arrow[r, "f"]
          & (\C^{n+1},0) \arrow[u, "j"]
          \end{tikzcd}
      \end{equation*}
commutes, where $i$ and $j$ are the natural inclusions. This diagram induces another commutative diagram
\begin{equation*}
        \begin{tikzcd}
          \O_{n+r}  \arrow[d, "i^*"]
          & \O_{n+1+r} \arrow[l, "F^*"]\arrow[d, "j^*"] \\
          \O_{n} 
          & \O_{n+1}.\arrow[l, "f^*"]
          \end{tikzcd}
      \end{equation*}

We seek to check that the above construction behaves properly under deformations. For the conductor ideal, one has that:
\begin{lemma}\label{lemma:conductor_specialises} If $C(F)$ is the conductor ideal of $F$, then $i^*(C(F))=C(f)$. 
\end{lemma}
\begin{proof} Let $(Z,0)$ denote the image of $F$. By definition, $C(F)=\{A\in \O_{Z,0} : A\cdot \O_{n+r} \subset \O_{Z,0}\}$. Take an element $A\in C(F)$, and let $b\in \O_n$. Since $i^*$ is surjective, it follows that $b=i^*(B)$ for some $B\in \O_{n+r}$. Hence, $i^*(A)b=i^*(AB)\in i^*(\O_{Z,0})=\O_{Y,0}$. Thus, $i^*(A)\in C(f)$. For the reverse inclusion, let $a\in C(f)$, and let $A\in \O_{n+r}$ such that $i^*(A)=a$. Then, for every $B\in \O_{n+r}$ one has that $AB\in (i^*)^{-1}(ai^*(B))\in (i^*)^{-1}(\O_{Y,0})\subset \O_{Z,0}$ due to the fact that $a\in C(f)$. The claim follows. \end{proof}

This means that the conductor ideal specialises properly under deformations. Moreover, it can be checked that $j^*(\F_1(F))=\F_1(f)$. Indeed, 
\begin{align*}
    j^*(\F_1(f))&=j^*((F^*)^{-1}(C(F))) = (f^*)^{-1} (i^*(C(F)))=\\ &=(f^*)^{-1}(C(f))=\F_1(f). 
\end{align*}

Let us denote $G\in \O_{n+1+r}$ to be a reduced equation of the image of $F$, in such a way that $j^*(G)=g$. Notice that $j^*(J(G))\neq J(g)$, due to the fact that $J(G)$ contains derivatives with respect to the parameters. If we consider the relative jacobian module $J_y(G)$ generated by the derivatives $\partial G / \partial y_i$ with respect to the variables $(y_1, \ldots, y_{n+1})\in \C^{n+1}$, then it is straightforward that $j^*(J_y(G))=J(g)$. 
\begin{definition} The map $F^*$ induces an epimorphism of $\O_{n+1+r}$-modules 
\begin{equation*}
    \dfrac{\F_1(F)}{J_y(G)}\rightarrow \dfrac{C(F)}{J_y(G)\cdot \O_{n+r}}.
\end{equation*}
We define $M_\rel(G)$ as the $\O_{n+1+r}$-module given as the kernel of the previous morphism. 
\end{definition}
\begin{remark} As in the previous case, the module $M_\rel(G)$ is determined by the exact sequence 
\begin{equation*}
    0\rightarrow M_\rel(G) \rightarrow \dfrac{\F_1(F)}{J_y(G)}\rightarrow \dfrac{C(F)}{J_y(G)\cdot \O_{n+r}} \rightarrow 0.
\end{equation*}
\end{remark}
\begin{remark}\label{remark:formulaMrelg} In practice, the formula of proposition \ref{prop:formulaMg} can be generalised to this relative module to yield
\begin{equation*}
    M_\rel (G)= \dfrac{(F^*)^{-1}(J(G)\cdot \O_{n+r})}{J_y (G)}.
\end{equation*}
The same proof works in this relative version, just by taking into account that $J_y (G)\cdot \O_{n+r} = J(G)\cdot \O_{n+r}$ (see lemma \ref{lemma:J=J_y}).  
\end{remark}
The goal of this section is to show that, if $n\geq 2$, the module $M_\rel(G)$ specialises to $M(g)$. In order to do so, we first consider $M_\rel (G)$ to be an $\O_r$-module through the natural inclusion $\O_r\subset \O_{n+1+r}$. With this, the aim is to show that 
\begin{equation*}
    M_\rel(G) \otimes \dfrac{\O_r}{\m_r} \cong M(g),
\end{equation*}
where $\cong$ denotes isomorphism as $\C$-vector spaces. Later in this section, we improve these results to perform a slightly different specialisation to obtain an isomorphism of $\O_{n+1}$-modules. \\

In order to prove that the given module specialises properly, we need to show a few lemmas first. 
\begin{lemma}\label{lemma:F1_specialises}
    If $F$ is an $r$-parameter unfolding of $f$, then 
    \begin{equation*}
        \F_1(F)\otimes \dfrac{\O_r}{\m_r} \cong \F_1(f).
    \end{equation*}
    Moreover, if $I=\m_r\cdot \O_{n+1+r}$, then $I\cdot \F_1(F) = I\cap \F_1(F)$. 
\end{lemma}
\begin{proof}
    Notice that $I$ is the kernel of the surjective map $j^*$. Thus, the first and second isomorphism theorems yield
    \begin{equation*}
        \dfrac{\F_1(F)+I}{I}\cong \dfrac{\F_1(F)}{\F_1(F)\cap I}\cong j^*(\F_1(F))=\F_1(f).
    \end{equation*}
    Hence,
    \begin{equation*}
        \dfrac{\O_{n+1+r}}{\F_1(F)}\otimes \dfrac{\O_r}{\m_r} = \dfrac{\O_{n+1+r}}{\F_1(f)+I}\cong \dfrac{\O_{n+1+r}/I}{(\F_1(F)+I)/I}\cong \dfrac{\O_{n+1}}{\F_1(f)}. 
    \end{equation*}
    Now, take the exact sequence of $\O_r$-modules 
    \begin{equation*}
        0\rightarrow \F_1(F)\rightarrow \O_{n+1+r}\rightarrow \dfrac{\O_{n+1+r}}{\F_1(F)}\rightarrow 0.
    \end{equation*}
    Notice that $\O_{n+1+r}/\F_1(F)$ is a determinantal ring of dimension $n+r-1$ (see remark \ref{remark:F1_determinantal}), whose fibre $\O_{n+1}/\F_1(f)$ has dimension $n-1$. Hence, $\O_{n+1+r}/\F_1(F)$ is an $\O_r$-flat module. This implies that the previous exact sequence can be tensored out to yield an exact sequence 
    \begin{equation*}
        0\rightarrow \F_1(F)\otimes \dfrac{\O_r}{\m_r}\rightarrow \O_{n+1} \rightarrow \dfrac{\O_{n+1}}{\F_1(f)}\rightarrow 0.
    \end{equation*}
    Hence, $\F_1(F)\otimes \dfrac{\O_r}{\m_r}\cong \F_1(f)$. For the second part, notice that 
    \begin{equation*}
        \F_1(f) \cong \dfrac{\F_1(F)+I}{I}\cong \dfrac{\F_1(F)}{\F_1(F)\cap I}. 
    \end{equation*}
    On the other hand, 
    \begin{equation*}
        \F_1(F)\otimes \dfrac{\O_r}{\m_r} \cong  \dfrac{\F_1(F)}{\m_r \cdot \F_1(F)} \cong \dfrac{\F_1(F)}{I\cdot \F_1(F)}. 
    \end{equation*}
    This forces $\F_1(F)\cap I=I\cdot \F_1(F)$ since both modules are isomorphic. 
\end{proof}
\begin{lemma}\label{lemma:C_specialises}
    If $F$ is an $r$-parameter unfolding of $f$, then 
    \begin{equation*}
        C(F)\otimes \dfrac{\O_r}{\m_r} \cong C(f). 
    \end{equation*}
    Moreover, if $J=\m_r\cdot \O_{n+r}$, then $J\cdot C(F) = J\cap C(F)$. 
\end{lemma}
\begin{proof}
    The proof is completely analogous to the one given in the previous case. The only place where both proofs differ is in the justification that $\O_{n+r}/C(F)$ is $\O_r$-flat. In this case, since $C(F)=(\Lambda)$ is a principal ideal for which $\O_{n+r}/C(F)$ has dimension $n+r-1$ with $n-1$ dimensional fibre $\O_{n}/C(f)$. The $\O_r$-flatness of $\O_{n+r}/C(F)$ then follows. 
\end{proof}
\begin{proposition}\label{prop:F/J_specialises} For any $r$-parameter unfolding $F$ of $f$, 
\begin{equation*}
    \dfrac{\F_1(F)}{J_y(G)}\otimes  \dfrac{\O_r}{\m_r} \cong \dfrac{\F_1(f)}{J(g)} \quad\text{and}\quad \dfrac{C(F)}{J_y(G)\cdot \O_{n+r}}\otimes  \dfrac{\O_r}{\m_r} \cong \dfrac{C(f)}{J(g)\cdot \O_{n}}.
\end{equation*}
\end{proposition}
\begin{proof} In order to simplify the notation, let us write $\F=\F_1(F)$ and $J=J_y(G)$. From lemma \ref{lemma:F1_specialises} and a combination of the isomorphism theorems it follows that 
\begin{align*}
    \dfrac{\F}{J}\otimes \dfrac{\O_r}{\m_r} &\cong \dfrac{\F/J}{I\cdot (\F/J)} = \dfrac{\F/J}{(I\cdot \F + J)/J} \cong \dfrac{\F}{I\cap \F + J} \cong \\
    &\cong \dfrac{\F/ (I\cap \F)}{(I\cap \F+J)/(I\cap \F)} = \dfrac{\F /(I\cap \F)}{J/(I\cap \F\cap J)} = \\
    &= \dfrac{\F/(I\cap \F)}{J/(I\cap J)}\cong \dfrac{(\F+I)/I}{(J+I)/J}\cong \dfrac{\F_1(f)}{J(g)}. 
\end{align*}
This shows the first part of the proposition. The second can be shown with a completely analogous process. 
\end{proof}
The following lemma is a key result to show that the module specialises properly. A few words of caution are in order here, since this result no longer holds when the domain of $f$ is non-smooth, as it will be commented in the following section. 
\begin{lemma}\label{lemma:J=J_y}
    For any $r$-parameter unfolding $F$ of $f$, $J_y (G)\cdot \O_{n+r} = J(G) \cdot \O_{n+r}$. 
\end{lemma}
\begin{proof} Denote by $F(x,u)=(f_u(x),u)$. In particular, the jacobian matrix of $F$ has the form 
\begin{equation*}
    dF=\begin{pmatrix} df_u & * \\ 0 & I_r \end{pmatrix}.
\end{equation*}
Denote by $M_1, \ldots, M_{n+1}, M'_1, \ldots, M'_r$ to the $n+r$-minors of $df$ so that $M_1, \ldots, M_{n+1}$ are the $n$-minors of $df_u$. It then follows that $M'_1, \ldots, M'_r$ can be generated from the minors $M_1, \ldots, M_{n+1}$, that is, there exist $a_{k,l}\in \O_{n+r}$ such that 
\begin{equation*}
    M'_k=\sum_{l=1}^{n+1} a_{kl}M_l
\end{equation*}
for each $k\in \{1,\ldots, r\}$. By Piene's lemma \ref{lemma:piene}, it follows that, up to the sign, 
\begin{equation*}
    \dfrac{\partial G}{\partial y_l}\circ F = \Lambda M_l, \quad
    \dfrac{\partial G}{\partial u_k}\circ F = \Lambda M'_k, 
\end{equation*}
where $\Lambda$ is the generator of $C(F)$. Thus, 
\begin{equation*}
    \dfrac{\partial G}{\partial u_k}\circ F = \sum_l a_{kl} \,\dfrac{\partial G}{\partial y_l}\circ F,
\end{equation*}
forcing $J(G)\cdot \O_{n+r} = J_y (G) \cdot \O_{n+r}$. 
\end{proof}
\begin{lemma}\label{lemma:C/J_CM} If $n\geq 2$, then $C(F)/(J_y(G)\cdot \O_{n+r})$ is Cohen-Macaulay of dimension $n+r-2$.
\end{lemma}
\begin{proof} Denote $C(F)=(\Lambda)$. Then, multiplication by $\Lambda$ yields an isomorphism $\O_{n+r}\rightarrow C(F)$ that maps $R(F)$ to $J(G)\cdot \O_{n+r}$ in virtue of Piene's lemma. Therefore, 
\begin{equation*}
    \dfrac{C(F)}{J_y (G)\cdot \O_{n+r}} \cong \dfrac{\O_{n+r}}{R(F)}.
\end{equation*}
Since $F$ is finite and $n\geq 2$, then $\O_{n+r}/R(F)$ is a determinantal ring of dimension $n+r-2$ (see \cite{juanjo}), and hence a Cohen-Macaulay module. 
\end{proof}
With all there lemmas being shown, we are now able to show the main result of the section, namely that the module $M_\rel (G)$ specialises to $M(g)$ provided $n\geq 2$. 
\begin{theorem} If $n\geq 2$, the module $M_\rel(G)$ specialises to $M(g)$. More formally,
\begin{equation*}
    M_\rel(G)\otimes \dfrac{\O_r}{\m_r}\cong M(g).
\end{equation*}
\end{theorem}
\begin{proof} From the definition of $M_\rel (G)$ we have the exact sequence 
\begin{equation*}
    0\rightarrow M_\rel(G) \rightarrow \dfrac{\F_1(F)}{J_y(G)}\rightarrow \dfrac{C(F)}{J_y(G)\cdot \O_{n+r}} \rightarrow 0.
\end{equation*}
Since $C(F)/(J_y(G)\cdot \O_{n+r})$ is a Cohen-Macaulay ring of dimension $n+r-2$ with fibre $C(f)/(J(g)\cdot \O_{n})$ of dimension $n-2$, it follows that $C(F)/(J_y(G)\cdot \O_{n+r})$ is $\O_r$-flat. Hence, tensoring with $\O_r/\m_r$ in the previous exact sequence of $\O_r$-modules yields the exact sequence
\begin{equation*}
    0\rightarrow M_\rel(G)\otimes \dfrac{\O_r}{\m_r} \rightarrow \dfrac{\F_1(f)}{J(g)}\rightarrow \dfrac{C(f)}{J(g)\cdot \O_{n}} \rightarrow 0.
\end{equation*}
It then follows that 
\begin{equation*}
    M_\rel(G)\otimes \dfrac{\O_r}{\m_r} \cong M(g),
\end{equation*}
since the above exact sequence fully determines $M(g)$ up to isomorphism. 
\end{proof}

Our last aim is to slightly improve the specialisation process to obtain an isomorphism of $\O_{n+1}$-modules. If we restrict $M_\rel (G)$ to be an $\O_r$-module and we perform the tensor product $\O_r/\m_r$, we naturally obtain a $\C$-vector space that is isomorphic to $M(g)$. However, this specialisation process ignores the fact that $M(g)$ is an $\O_{n+1}$-module. In order to take this fact into account, one should perform the specialisation process in the following way: since $M_\rel (G)$ is an $\O_{n+1+r}$-module, the tensor product 
\begin{equation*}
    M_\rel(G)\otimes \dfrac{\O_{n+1+r}}{\m_r\cdot \O_{n+1+r}}
\end{equation*}
provides naturally an $\O_{n+1}$-module due to the fact that $\O_{n+1+r}/\m_r\cdot \O_{n+1+r}\cong \O_{n+1}$, where $\m_r=(u_1, \ldots , u_r)$ is the maximal ideal of $\O_r$, generated by the parameters of the unfolding. In addition, one has that 
\begin{equation*}
    M_\rel(G)\otimes \dfrac{\O_{n+1+r}}{\m_r\cdot \O_{n+1+r}} \cong \dfrac{M_\rel(G)}{\m_r\cdot M_\rel (G)}
\end{equation*}
embodies the spirit of forcing the parameters of the unfolding to be equal to 0. Hence, this analogous process of specialisation provides a similar result, but taking the $\O_{n+1}$-module structure into account. \\

In what follows, we show that 
\begin{equation*}
    M_\rel(G)\otimes \dfrac{\O_{n+1+r}}{\m_r\cdot \O_{n+1+r}} \cong M(g),
\end{equation*}
where $\cong $ now denotes isomorphism of $\O_{n+1}$-modules. \\

In order to prove it, we revisit each of the required lemmas and propositions that we have applied in the previous case to obtain this result:

\begin{lemma} If $F$ is an unfolding of $f$, then 
$$\F_1(F)\otimes \dfrac{\O_{n+1+r}}{\m_r\cdot \O_{n+1+r}} \cong \F_1(f)$$
as $\O_{n+1}$-modules. 
\end{lemma}
\begin{proof}
    The same proof that was given in lemmma \ref{lemma:F1_specialises} shows that 
    \begin{equation*}
        \dfrac{\O_{n+1+r}}{\F_1(F)}\otimes \dfrac{\O_{n+1+r}}{\m_r\cdot \O_{n+1+r}} \cong \dfrac{\O_{n+1}}{\F_1(f)}.
    \end{equation*}
    Moreover, since $I\cap \F_1(F)= I\cdot \F_1(F)$, it follows by proposition \ref{prop:torsionIJ} that 
    \begin{equation*}
        \Tor_1^{\O_{n+1+r}} \left( \dfrac{\O_{n+1+r}}{\F_1(F)}, \dfrac{\O_{n+1+r}}{\m_r\cdot \O_{n+1+r}} \right) = \dfrac{I\cap \F_1(F)}{I\cdot \F_1(F)}=0.
    \end{equation*}
    Hence, the short exact sequence of $\O_{n+1+r}$-modules
    \begin{equation*}
        0 \rightarrow \F_1(F)\rightarrow \O_{n+1+r}\rightarrow \dfrac{\O_{n+1+r}}{\F_1(F)}\rightarrow 0
    \end{equation*}
    yields after tensoring that the sequence
    \begin{equation*}
        0 \rightarrow \F_1(F)\otimes \dfrac{\O_{n+1+r}}{\m_r\cdot \O_{n+1+r}}\rightarrow \O_{n+1}\rightarrow \dfrac{\O_{n+1}}{\F_1(f)}\rightarrow 0
    \end{equation*}   
    is short exact as well. Therefore, 
$$\F_1(F)\otimes \dfrac{\O_{n+1+r}}{\m_r\cdot \O_{n+1+r}} \cong \F_1(f)$$
    as $\O_{n+1}$-modules. 
\end{proof}
An analogous proof immediately shows the same for the conductor ideal, with the only consideration that $C(f)$ is naturally an $\O_n$-module. 
\begin{lemma} If $F$ is an unfolding of $f$, then 
$$C(F)\otimes \dfrac{\O_{n+r}}{\m_r\cdot \O_{n+r}} \cong C(f)$$
as $\O_{n}$-modules. 
\end{lemma}
Furthermore, by applying these two lemmas and performing the same reasoning as was done in proposition \ref{prop:F/J_specialises} we deduce that:
\begin{proposition} For any $r$-parameter unfolding $F$ of $f$, 
\begin{equation*}
    \dfrac{\F_1(F)}{J_y(G)}\otimes  \dfrac{\O_{n+1+r}}{\m_r\cdot \O_{n+1+r}} \cong \dfrac{\F_1(f)}{J(g)}
\end{equation*}
and 
\begin{equation*}
    \dfrac{C(F)}{J_y(G)\cdot \O_{n+r}}\otimes  \dfrac{\O_{n+r}}{\m_r\cdot \O_{n+r}} \cong \dfrac{C(f)}{J(g)\cdot \O_{n}}.
\end{equation*}
\end{proposition}

With these results, we are now able to show that the module $M_\rel (G)$ specialises properly when the $\O_{n+1}$-module structure is preserved. 
\begin{theorem} If $n\geq 2$, then
\begin{equation*}
    M_\rel(G)\otimes \dfrac{\O_{n+1+r}}{\m_r\cdot \O_{n+1+r}} \cong M(g),
\end{equation*}
as $\O_{n+1}$-modules. 
\end{theorem}
\begin{proof}
    Notice first that the same proof that was given to show that $\F_1(F)$ specialises properly in the two different meanings that have been studied here can be applied to the ramification ideal $R(F)$. Indeed, $j^*(R(F))=R(f)$, and $\O_{n+1+r}/R(F)$ is a determinantal ring provided $n\geq 2$ as was mentioned before. Then, the same reasonings provided for $\F_1(F)$ show that 
    $$R(F)\otimes \dfrac{\O_{n+1+r}}{\m_r\cdot \O_{n+1+r}} \cong R(f)$$
    as $\O_{n+1}$-modules, and, in particular, that if $I=\m_r\cdot \O_{n+1+r}$, then $I\cdot R(F) = I\cap R(F)$. Hence, 
    \begin{multline*}
        \Tor_1^{\O_{n+1+r}} \left( \dfrac{C(F)}{J_y(G)\cdot \O_{n+r}}, \dfrac{\O_{n+1+r}}{\m_r\cdot \O_{n+1+r}} \right) \cong \\ \cong  \Tor_1^{\O_{n+1+r}} \left( \dfrac{\O_{n+1+r}}{R(F)}, \dfrac{\O_{n+1+r}}{\m_r\cdot \O_{n+1+r}} \right) = \dfrac{I\cap R(F)}{I\cdot R(F)}=0.
    \end{multline*}
    Therefore, from the short exact sequence of $\O_{n+1+r}$-modules that defines $M_\rel (G)$ by
    \begin{equation*}
    0\rightarrow M_\rel(G) \rightarrow \dfrac{\F_1(F)}{J_y(G)}\rightarrow \dfrac{C(F)}{J_y(G)\cdot \O_{n+r}} \rightarrow 0,
    \end{equation*}
    one obtains after tensoring that the sequence
        \begin{equation*}
    0\rightarrow M_\rel(G) \otimes \dfrac{\O_{n+1+r}}{\m_r\cdot \O_{n+1+r}} \rightarrow \dfrac{\F_1(f)}{J(g)}\rightarrow \dfrac{C(f)}{J(g)\cdot \O_{n}} \rightarrow 0
    \end{equation*}
    is short exact. Thus, 
    \begin{equation*}
    M_\rel(G)\otimes \dfrac{\O_{n+1+r}}{\m_r\cdot \O_{n+1+r}} \cong M(g),
    \end{equation*}
    as $\O_{n+1}$-modules, and the result follows. 
\end{proof}

\newpage
\section{A Jacobian module for distentanglements of mappings defined on \icis}

In this chapter, we generalise the definition that has been introduced in the previous chapter for mappings defined on \icis with the intention to show the Mond conjecture in this general framework for $n=2$. We then follow the ideas of \cite{BobadillaNunoPenafort} to extend them, and we prove the apply them to prove the desired case of the conjecture. 

\subsection{Definition of the generalised version of the module \textit{M(g)}}
In this section, a general version of the module is defined for mappings with non-smooth source. A few words of caution are in order here, since some results that hold in the smooth case are no longer true in this setting. We will comment on them, and we provide a natural way to define this general version of the module in terms of the construction developed in the smooth source case. \\

Let $f:(X,S)\rightarrow (\C^{n+1},0)$ be an $\A$-finite mapping defined on an \icis $(X,S)\subset (\C^{n+k},S)$ of dimension $n$. As $f$ is finite, its image is a hypersurface of $(\C^{n+1},0)$, and hence it can be described by a reduced equation $g\in \O_{n+1}$. Write $(X,S)=h^{-1}(0)$, where $h:(\C^{n+k},S)\rightarrow (\C^k,0)$ is a $\K$-finite mapping. Consider an analytic extension $\tilde{f}:(\C^{n+k},S)\rightarrow (\C^{n+1},0)$ of $f$, and write 
$$\hat{f}=(\Tilde{f},h):(\C^{n+k},S)\rightarrow (\C^{n+k+1},0).$$ It is therefore clear that the restriction of $\hat{f}$ to $(X,S)$ is precisely the mapping $(f,0)$. Hence, $\hat{f}$ is a finite mapping, since $\hat{f}^{-1}(0)=f^{-1}(0)=S$. Moreover, the diagram
\begin{center}
    \begin{tikzcd}
	& {(\C^k,0)} \\
	{(\C^{n+k},0)} & {(\C^{n+k+1},0)} \\
	{(X,S)} & {(\C^{n+1},0)}
	\arrow["f", from=3-1, to=3-2]
	\arrow["i", hook, from=3-1, to=2-1]
	\arrow["\hat{f}", from=2-1, to=2-2]
	\arrow["j"', hook, from=3-2, to=2-2]
	\arrow["h", two heads, from=2-1, to=1-2]
	\arrow[two heads, from=2-2, to=1-2]
\end{tikzcd}
\end{center}
commutes, where $i$ is the inclusion and $j$ is the obvious immersion. In particular, $\hat{f}$ is an unfolding of $f$ deforming both the mapping and the domain. After taking representatives, the induced deformations of $f$ are the mappings $\hat{f}_t:X_t\subset \C^{n+k}\rightarrow \C^{n+1}$ defined as $\hat{f}_t(x)=\tilde{f}(x,t)$, where $X_t=h^{-1}(t)$ for $t\in \C^k$ small enough. \\

The key idea is that $\hat{f}$ is, in some sense, the simplest unfolding of $f$ with smooth source. Hence, the definition of the module for the mapping $f$ will be performed through a specialisation of the one from $\hat{f}$. \\

Since $\hat{f}$ is a finite mapping, its image $\im \hat{f}$ is a hypersurface of $(\C^{n+k+1},0)$. Let us consider $\hat{g} \in \O_{n+k+1}$ to be a reduced equation of the image such that $\hat{g} \circ j =g$. As we have already mentioned in the previous section, a result from Piene shows that $J(\hat{g})\cdot \O_{n+k}\subset C(\hat{f})$, where $J(\hat{g})\cdot \O_{n+k} = \hat{f}^*(J(\hat{g}))$. On the other hand, since $\hat{f}$ is a finite mapping with degree 1 onto its image, a result from Mond-Pellikaan shows that $\F_1(\hat{f}) \cdot \O_{n+k} = C(\hat{f})$. \\

Let us denote by $(y,z)$ to the coordinates in $(\C^{n+k+1},0)$, with $y\in \C^{n+1}$ and $z\in \C^k$, and consider the jacobian ideal 
\begin{equation*}
    J_y(\hat{g} ) = \Big\langle \dfrac{\partial \hat{g}}{\partial y_1}, \ldots, \dfrac{\partial \hat{g}}{\partial y_{n+1}}\Big\rangle
\end{equation*}
generated by the derivatives with respect to $y=(y_1, \ldots, y_{n+1})$. In this case, it follows that the restriction of $\hat{f}^*$ to $\F_1(\hat{f})$ induces an epimorphism of $\O_{n+k+1}$-modules 
\begin{equation*}
    \dfrac{\F_1(\hat{f})}{J_y(\hat{g})}\rightarrow \dfrac{C(\hat{f})}{J(\hat{g})\cdot \O_{n+k}}.
\end{equation*}
We define $N(\hat{g})$ to be the $\O_{n+k+1}$-module given by the kernel of this morphism, and define 
\begin{equation*}
    M(g)=N(\hat{g})\otimes\dfrac{\O_{n+k+1}}{\m_k\cdot \O_{n+k+1}},
\end{equation*}
which has an $\O_{n+1}$-module structure. Hence, $M(g)$ is defined by taking into account that $\hat{f}$ is an unfolding of $f$.

\begin{remark} Notice that the module $N(\hat{g})$ is not exactly the same as the module of the previous section. Indeed, the source of the morphism that defines $N(\hat{g})$ has the jacobian ideal $J_y(\hat{g})$ where only partial derivatives with respect to the parameters of $(Y,0)$ are taken into consideration, while the module in the target does have all the partial derivatives. Altough this may seem to be somewhat whimsical, this will be shown to be exactly what is needed for the module to specialise properly. 
\end{remark}
Notice that this module $M(g)$ coincides with the given in the article of Bobadilla-Nuño-Peñafort in the smooth case just by taking $\hat{f}=f$  and $k=0$. \\
\begin{remark} It is crucial to notice that, in general, $J(\hat{g})\cdot \O_{n+k}\neq J_y(\hat{g})\cdot \O_{n+k}$, in contrast with what happens in the smooth case. This important fact is what infuences this definition for $M(g)$. 
\end{remark}
\begin{remark} The module $N(\hat{g})$ is determined by the short exact sequence
\begin{equation*}
    0\rightarrow N(\hat{g}) \rightarrow \dfrac{\F_1(\hat{f})}{J_y(\hat{g})}\rightarrow \dfrac{C(\hat{f})}{J(\hat{g})\cdot \O_{n+k}}\rightarrow 0.
\end{equation*}
After tensoring with $\O_{n+k+1}/\m_k\cdot \O_{n+k+1}$, it follows that the module $M(g)$ fits in the exact sequence 
\begin{equation*}
    M(g) \rightarrow \dfrac{\F_1(f)}{J(g)}\rightarrow \dfrac{C(\hat{f})}{J(\hat{g})\cdot \O_{n+k}}\otimes\dfrac{\O_{n+k+1}}{\m_k\cdot \O_{n+k+1}}\rightarrow 0.
\end{equation*}
Notice that 
\begin{equation*}
    \dfrac{C(\hat{f})}{J(\hat{g})\cdot \O_{n+k}} \cong \dfrac{\O_{n+k}}{R\hat{f}},
\end{equation*}
where $R\hat{f}$ is the ramification ideal of $\hat{f}$, which is a determinantal ring of dimension $n+k-2$. Hence, 
\begin{equation*}
    \dfrac{C(\hat{f})}{J(\hat{g})\cdot \O_{n+k}}\otimes\dfrac{\O_{n+k}}{\m_k\cdot \O_{n+k}} \cong \dfrac{\O_{X}}{i^*(R\hat{f})}\cong \dfrac{\O_{n+k}}{R\hat{f}+(h_1, \ldots, h_k)}.
\end{equation*}
Notice that, outside the origin, $f$ is stable. Then, the previous module coincides with the ramification ideal of a stable mapping, which has dimension $n-2$. Hence, the module has dimension $n-2$, as the addition of the origin does not perturb the dimension. In particular, the module $C(\hat{f})/J(\hat{g})\cdot \O_{n+k}$ is $\O_k$-flat, which forces the sequence 
\begin{equation*}
    0\rightarrow M(g) \rightarrow \dfrac{\F_1(f)}{J(g)}\rightarrow \dfrac{\O_{X}}{i^*(R\hat{f})}\rightarrow 0
\end{equation*}
to be short exact. Hence, $M(g)$ is also determined by a short exact sequence, and this could also be taken as a definition of $M(g)$ in this case. 
\end{remark}
\begin{proposition} The following formula holds:
\begin{equation*}
    N(\hat{g})=\dfrac{(\hat{f}^*)^{-1} (J(\hat{g})\cdot \O_{n+k})}{J_y(\hat{g})}. 
\end{equation*}
\end{proposition}
\begin{proof} By definition, it is clear that 
\begin{equation*}
    N(\hat{g})=\dfrac{(\hat{f}^*)^{-1} (J(\hat{g})\cdot \O_{n+k})\cap \F_1(\hat{f})}{J_y(\hat{g})}.
\end{equation*}
Since $J(\hat{g})\cdot \O_{n+k}\subset C(\hat{f})$, then $(\hat{f}^*)^{-1} (J(\hat{g})\cdot \O_{n+k})\subset  \F_1(\hat{f})$, and the first formula follows.
\end{proof}

\subsection{A relative version for the generalised module}

In the following section, a relative version of the module for unfoldings is defined. For the sake of simplicity, we consider unfoldings $F:(\C^{n+r+k},S\times 0)\rightarrow (\C^{n+1+r+k},0)$ of the mapping $\hat{f}$ instead of general unfoldings of $f$ with possibly non-smooth source. Let $G\in \O_{n+1+r+k}$ be an equation for the image of $F$ such that $G(y,z,0)=\hat{g}(y,z)$, where $(y,z,u)$ are the coordinates of $(\C^{n+1+r+k},0)$, with $y\in \C^{n+1}, z\in \C^k$ and $u\in \C^r$. We then have that the diagram
\[\begin{tikzcd}
	{(\C^{n+r+k},S\times 0)} & {(\C^{n+1+r+k},0)} \\
	{(\C^{n+k},S)} & {(\C^{n+k+1},0)} & {(\C,0)} \\
	{(X,S)} & {(\C^{n+1},0)}
	\arrow["{\hat{f}}", from=2-1, to=2-2]
	\arrow["f", from=3-1, to=3-2]
	\arrow["i", hook, from=3-1, to=2-1]
	\arrow["j"', hook, from=3-2, to=2-2]
	\arrow["{\hat{i}}", hook, from=2-1, to=1-1]
	\arrow["{\hat{j}}"', hook, from=2-2, to=1-2]
	\arrow["F", from=1-1, to=1-2]
	\arrow["g", from=3-2, to=2-3]
	\arrow["{\hat{g}}", from=2-2, to=2-3]
	\arrow["G", from=1-2, to=2-3]
\end{tikzcd}\]
commutes, where $\hat{i}, \hat{j}$ are the obvious immersions. Let us denote by 
\begin{equation*}
    J_y(G ) = \Big\langle \dfrac{\partial G}{\partial y_1}, \ldots, \dfrac{\partial G}{\partial y_{n+1}}\Big\rangle, \,\,
    J_z(G ) = \Big\langle \dfrac{\partial G}{\partial z_1}, \ldots, \dfrac{\partial G}{\partial z_k}\Big\rangle, 
\end{equation*}
and $J_{y,z}(G)=J_y(G)+J_z(G)$. We therefore define $M_\rel (G)$ as the kernel of the module epimorphism 
\begin{equation*}
    \dfrac{\F_1(F)}{J_y (G)}\rightarrow \dfrac{C(F)}{J_{y,z}(G)\cdot \O_{n+r+k}}.
\end{equation*}
In this case, the given definition of the $\O_{n+1+r+k}$-module $M_\rel (G)$ properly specialises to $M(g)$ seen as an $\O_{n+1}$-module, in the sense of the following theorem:
\begin{theorem}
    If $n\geq 2$, then 
    \begin{equation*}
        M_\rel (G) \otimes \dfrac{\O_{n+1+r+k}}{\m_{r+k}\cdot \O_{n+1+r+k}}\cong M(g).
    \end{equation*}
\end{theorem}
\begin{proof} The same proof of the smooth case shows that 
\begin{equation*}
    M_\rel (G)\otimes \dfrac{\O_{n+1+r+k}}{\m_r\cdot \O_{n+1+r+k}}\cong N(\hat{g}).
\end{equation*}
Indeed, the only difference is that the definition of $N(\hat{g})$ has a term $J_y(\hat{g})$, but it holds that $\hat{j}(J_y(G))=J_y(\hat{g})$ and, as in the smooth case, $\hat{j}(J_{y,z}(G))=J(\hat{g})$. Hence, 
\begin{multline*}
   M_\rel (G) \otimes \dfrac{\O_{n+1+r+k}}{\m_{r+k}\cdot \O_{n+1+r+k}} = \\ = \left( M_\rel (G) \otimes \dfrac{\O_{n+1+r+k}}{\m_{r}\cdot \O_{n+1+r+k}}\right)\otimes \dfrac{\O_{n+1+k}}{\m_{k}\cdot \O_{n+1+k}} = \\ = N(\hat{g})\otimes \dfrac{\O_{n+1+k}}{\m_{k}\cdot \O_{n+1+k}}  =M(g),
\end{multline*}
where the last equality holds by definition of $M(g)$. 
\end{proof}
Despite what happened with the module $M(g)$, the relative version $M_\rel (G)$ is fully determined as the kernel of an epimorphism. Hence, it fits into the short exact sequence 
\begin{equation*}
    0\rightarrow M_\rel (G) \rightarrow \dfrac{\F_1(F)}{J_y (G)}\rightarrow \dfrac{C(F)}{J_{y,z}(G)\cdot \O_{n+r+k}}\rightarrow 0.
\end{equation*}
Furthermore, one has that
\begin{proposition} The following formula holds:
\begin{equation*}
    M_\rel (G)=\dfrac{(F^*)^{-1} (J_{y,z}(G)\cdot \O_{n+r+k})}{J_y(G)}.
\end{equation*}
\end{proposition}
In contrast with the absolute case, we have that $J_{y,z}(G)\cdot \O_{n+r+k}=J(G)\cdot \O_{n+r+k}$. Hence, they can be interchanged indistinctively in this formula. \\


If one compares this case with what has been done in the smooth case, it is relevant to notice that the specialisation process has been performed in a different way. Indeed, in this case, the double specialization requires to treat work with these modules in a more natural way. When we restrict $M_\rel (G)$ to be an $\O_{r+k}$-module via the natural inclusion $\O_{r+k}\rightarrow \O_{n+1+r+k}$ induced by the projection $\C^{n+1}\times \C^{r+k}\rightarrow \C^{r+k}$, it turns out that the vector space structure of both specialisations coincide, as the following propostion asserts:
\begin{proposition}\label{prop:specialisation} If $M_\rel(G)$ is considered an $\O_{r+k}$-module, then 
\begin{equation*}
    M_\rel(G)\otimes \dfrac{\O_{r+k}}{\m_{r+k}}\cong M(g)
\end{equation*}
as $\C$-vector spaces. 
\end{proposition}
\begin{proof}
    Notice that 
    \begin{equation*}
    M_\rel(G)\otimes \dfrac{\O_{r+k}}{\m_{r+k}}\cong \dfrac{M_\rel (G)}{\m_{r+k}\cdot M_\rel (G)}
\end{equation*}
as $\O_{r+k}$-modules. On the other hand, it has already been checked in the previous theorem that  
\begin{equation*}
    M(g)\cong M_\rel (G) \otimes \dfrac{\O_{n+1+r+k}}{\m_{r+k}\cdot \O_{n+1+r+k}}
\end{equation*}
as $\O_{n+1}$-modules. Furthermore, 
\begin{equation*}
    M_\rel (G) \otimes \dfrac{\O_{n+1+r+k}}{\m_{r+k}\cdot \O_{n+1+r+k}} \cong \dfrac{M_\rel (G)}{(\m_{r+k}\cdot \O_{n+1+r+k})\cdot M_\rel (G)}
\end{equation*}
naturally as $\O_{n+1+r+k}$-modules. Notice that $$(\m_{r+k}\cdot \O_{n+1+r+k})\cdot M_\rel (G)=\m_{r+k}\cdot M_\rel (G).$$ Indeed, both sets are given by the elements of the form $\sum_i u_i m_i$, where $\m_{r+k}=(u_1, \ldots, u_{r+k})$, and for some $m_i\in M_\rel (G)$. Hence, 
\begin{equation*}
     \dfrac{M_\rel (G)}{\m_{r+k}\cdot M_\rel (G)} = \dfrac{M_\rel (G)}{(\m_{r+k}\cdot \O_{n+1+r+k})\cdot M_\rel (G)}.
\end{equation*}
Connecting all the previous isomorphisms and equalities shows that $M(g)$ is isomorphic as a $\C$-module to $M_\rel(G)\otimes (\O_{r+k}/\m_{r+k})$, and hence completes the proof.
\end{proof}
Therefore, the specialisation process can be performed either restricting first $M_\rel (G)$ to be an $\O_{r+k}$-module, or keeping its whole $\O_{n+1+r+k}$-module structure. The obtained modules are, respectively, $M_\rel (G)\otimes (\O_{n+1+r+k}/\m_{r+k}\cdot \O_{n+1+r+k})$ and $M_\rel (G)\otimes (\O_{r+k}/\m_{r+k})$. In both scenarios one recovers $M(g)$ in some sense: in the former case, the result is an $\O_{n+1}$-module isomorphic to $M(g)$, and, in the latter, a $\C$-module isomorphic to $M(g)$. Hence, the specialisation procces that one should perform depends on whether one needs to keep the module structure or not. For most of the cases, the only information that needs to be keeped is the complex dimension as a vector space, and hence both methods are valid. \\

Lastly, an important result regarding the form of the module $M_\rel (G)$ when $F$ is stable is the following:
\begin{proposition} Let $F$ be a stable unfolding of $\hat{f}$ and $G$ an equation such that $G\in J(G)$. Then, 
\begin{equation*}
    M_\rel (G)= \dfrac{J(G)}{J_y(G)}.
\end{equation*}
\end{proposition}
\begin{proof}
    Since $F$ is stable and $G\in J(G)$, the smooth version of the module $M(G)$ from the article of Bobadilla-Nuño-Peñafort satisfies that $M(G)=0$. In this case, the formula
    \begin{equation*}
        M(G)=\dfrac{(F^*)^{-1}(J(G)\cdot \O_{n+r+k})}{J(G)}
    \end{equation*}
    holds, and hence $(F^*)^{-1}(J(G)\cdot \O_{n+r+k})=J(G)$. Furthermore, in the smooth case it follows that $J(G)\cdot \O_{n+r+k}=J_{y,z}(G)\cdot \O_{n+r+k}$. Hence, formula of the previous proposition yields that
    \begin{equation*}
    M_\rel (G)=\dfrac{(F^*)^{-1} (J_{y,z}(G)\cdot \O_{n+r+k})}{J_y(G)}=\dfrac{J(G)}{J_y(G)},
\end{equation*}
and the claim follows. 
\end{proof}
\subsection{Relation between $\dim_\C M(g)$ and  $\codimAe (X,f)$}
Let $F$ be a stable unfolding of $\hat{f}$, and consider an equation $G$ of the image of $F$, namely, $(Z,0)$, that satisfies $G\in J(G)$. Then, the last result of the previous section showed that $M_\rel (G) =J(G)/J_y(G)$. Let us relate this with the codimension of $(X,f)$, which can be determined through the formula
\begin{equation*}
    \codimAe (X,f)=\dim_\C \dfrac{\theta(i)}{ti(\theta_{n+1})+i^*(\Derlog Z)},
\end{equation*}
where $Y$ denotes the image of $F$, and $i:(\C^{n+1},0)\rightarrow (\C^{n+1+k+r},0)$ denotes the natural immersion (see \cite{Mond-Montaldi} for more details). \\

Recall that $\Derlog Z = \{\xi \in \theta_{n+1+r+k}: \xi(G)=\lambda (G)\}$, and $\Derlog G = \{\xi \in \theta_{n+1+r+k}: \xi(G)=0\}$. Notice that $G\in J(G)$, so that $G=\sum_{s} a_s \partial_s G $, where $\partial_s G$ denotes the partial derivatives with respect to all the variables in $(y,z,u)\in \C^{n+1+r+k}$. Hence, the vector field $\epsilon = \sum_s a_s\partial_s$ satisfies that $\epsilon (G)=G$, where $\partial_s$ denotes the coordinate vector field associated with the $s$-th coordinate, where $s\in \{1, \ldots, n+1+r+k\}$. Furthermore, $\Derlog Z = \Derlog G \oplus \langle \epsilon \rangle$. Indeed, for the direct inclusion, notice that each $\xi$ with $\xi(G)=\lambda G$ can be rewritten via $\eta:=\xi-\lambda \epsilon$. Since $\eta (G)=\lambda G-\lambda G=0$, then $\eta \in \Derlog G$ and $\xi= \eta +\lambda \epsilon$. The reverse inclusion is obvious. Lastly, the sum is direct due to the fact that each vector field $\xi$ satisfying both $\xi=\lambda \epsilon$ and $\xi (G)=0$ is forced to have that $0=\xi(G)=\lambda G$, so that $\lambda =0$ and $\xi=0$. \\

We therefore have that  
\begin{equation*}
    \codimAe (X,f)=\dim_\C \dfrac{\theta(i)}{ti(\theta_{n+1})+i^*(\Derlog G)+i^*(\epsilon)}.
\end{equation*}
Notice that the evaluation mapping $\text{ev}:\theta_{n+1+r+k}\rightarrow J(G)$ given by $\xi \mapsto \xi (G)$ is a surjective mapping with kernel $\Derlog G$. Hence, it induces an isomorphism
\begin{equation*}
    \dfrac{\theta_{n+1+r+k}}{\Derlog G} \cong J(G).
\end{equation*}
Thus, 
\begin{equation*}
    \dfrac{\theta_{n+1+r+k}}{\langle \frac{\partial}{\partial y_1}, \ldots, \frac{\partial}{\partial y_{n+1}} \rangle +\Derlog G} \cong \dfrac{J(G)}{J_y (G)}=M_\rel (G).
\end{equation*}
Tensoring with $\O_{r+k}/\m_{r+k}$ yields that 
\begin{equation*}
    \dfrac{\theta(i)}{ti(\theta_{n+1}) +i^*\Derlog G} \cong M(g).
\end{equation*} 
Now, notice that the evaluation map acting on $\epsilon$ gives $\epsilon (G)=G$, and hence $i^*(\text{ev}(\epsilon))=i^*(G)=g$. Therefore, if $K(g)=(J(g)+(g))/J(g)$, then 
\begin{multline*}
    0 \rightarrow K(g) \rightarrow \dfrac{\theta(i)}{ti(\theta_{n+1}) +i^*\Derlog G} \rightarrow \\
    \rightarrow \dfrac{\theta(i)}{ti(\theta_{n+1}) +i^*\Derlog G +i^*(\epsilon)}\rightarrow 0
\end{multline*}
is a short exact sequence. Indeed, the evaluation map satisfies that $\text{ev}(ti(\theta_{n+1}))=J(g)$ and that $\text{ev}(i^*\Derlog G)=0$. Hence, the evaluation map yields an isomorphism 
\begin{equation*}
    \dfrac{ti(\theta_{n+1}) +i^*\Derlog G+i^*(\epsilon)}{ti(\theta_{n+1}) +i^*\Derlog G }\cong \dfrac{J(g)+(g)}{J(g)}=K(g).
\end{equation*}
After taking lengths in the exact sequence, and taking into account the previous assertions, it follows that 
\begin{equation*}
    \dim_\C M(g) = \dim_\C K(g)+\codimAe (X,f). 
\end{equation*}
This shows that $\dim_\C M(g)\geq \codimAe (X,f)$, with equality in case that $g$ is weighted homogeneous. Furthermore, this formula shows that $\dim_\C M(g)$ only depends on the isomorphism class of $g$, and neither depends on the mapping $f$ nor on the chosen extension $\hat{f}$. 

\subsection{The Jacobian module and the Mond conjecture}\label{section:5_3}
This section is the centerpiece of the project. We present one of the main results we aim to establish, namely, a formula for the image Milnor number expressed in terms of the Samuel multiplicity of the module $M_\rel(G)$. Additionally, we establish that the Mond conjecture holds  provided the module $M_\rel (G)$ is Cohen-Macaulay. In order to do so, we state a relevant result from Siersma \cite{Siersma} regarding the homotopy type of a fibre:
\begin{theorem}[Siersma] Let $g:(\C^{n+1},0)\rightarrow (\C,0)$ define a reduced hypersurface $(X_0,0)$, not necessarily with isolated singularity, and let $G:(\C^{n+1+r},0)\rightarrow (\C,0)$ be a deformation of $g$ such that
\begin{enumerate}
    \item $G$ is topologically trivial over the Milnor sphere $\partial B_\epsilon$, and
    \item for all $u$, all the critical points of $g_u$ which are not in $X_u=g_u^{-1}(0)\cap B_\epsilon$ are isolated. 
\end{enumerate}
    Then, $X_u\cap B_\epsilon$ is homotopy equivalent to a wedge of $n$-spheres and the number of such $n$-spheres is equal to 
    \begin{equation*}
        \sum_{y\in B_\epsilon \setminus X_u} \mu (g_u; y),
    \end{equation*}
    where $\mu(g_u,y)$ denotes the Milnor number of the function $g_u$ at the point $y$. 
\end{theorem}

To ensure that the reader is well-versed with the notion of Samuel multiplicity of a module, we have included the technical details in appendix \ref{appendix:2}. We highly recommend that the reader go through this appendix before proceeding to this section.
\begin{theorem}  If $F$ is either a stable unfolding or a stabilisation of $\hat{f}$ and with the notations made above, $$\mu_I(X,f)=e(\m_{r+k},M_\rel(G))$$
considering $M_\rel (G)$ an $\O_{r+k}$-module. 
\end{theorem}
\begin{proof} Take a representative of $F$ and let $w\in \C^k \times \C^r$ be a generic value. The conservation of multiplicity \ref{app:samuel1} implies that
    \begin{equation*}
        e\Big(\m_{r+k}, M_\rel (G)\Big)=\sum_{p\in B_\epsilon} e\Big( \m_{r+k,w}, M_\rel (G)_{(p,w)} \Big).
    \end{equation*}
    In order to compute the previous multiplicity, let us take first the points $p\in Y_w\cap B_\epsilon$, where $Y_w$ is the image of $f_w:X_w\rightarrow \C^{n+1}$. Since the module $M_\rel (G)$ specialises to $M(g)$, it follows that 
    \begin{equation*}
        M_\rel (G)_{(p,w)}\otimes \dfrac{\O_{r+k,w}}{\m_{r+k,w}}\cong M(g_w)_p
    \end{equation*}
    as $\C$-vector spaces in virtue of proposition \ref{prop:specialisation}.     Either if $F$ is a stable unfolding or a stabilisation, it follows that, provided $w$ is generic, $f_w$ is a stable mapping, and hence $ M(g_w)_p=0$ for each $p\in Y_w$. Hence, the points $p\in Y_w \cap B_\epsilon$ do not contribute to the term $e(\m_{r+k}, M_\rel (G))$. \\

    On the other hand, if $p\in B_\epsilon / Y_w$, then 
    \begin{equation*}
        M_\rel (G)_{(p,w)}  = \dfrac{\O_{\C^{k+r}\times B_\epsilon, (p,w)}}{J_y (G)}
    \end{equation*}
    is a module with dimension $\geq k+r$. Indeed, this follows from the exact sequence that defines $M_\rel (G)$, since the localisation of $C(F)/(J_{y,z}(G)\cdot \O_{n+k+r})$ is zero if $p\notin Y_w$, and $\F_1(F)$ localises to $\O_{\C^{k+r}\times B_\epsilon, (p,w)}$. Moreover, 
    \begin{equation*}
        M_\rel (G)_{(p,w)}\otimes \dfrac{\O_{r+k,w}}{\m_{r+k,w}} \cong  \dfrac{\O_{B_\epsilon, p}}{J(g_w)}
    \end{equation*}
    is a module with dimension $0$, since it has finite length due to the fact that $g_u$ has isolated singularity. This implies that the dimension of $M_\rel (G)$ is $\leq r+k$. Hence, $M_\rel (G)_{(p,u)}$ is a complete intersection ring, and in particular a Cohen-Macaulay $\O_{r+k}$-module of dimension $r+k$. Hence, 
    \begin{align*}
        e\Big( \m_{r+k,w}, M_\rel (G)_{(p,w)} \Big) &= \dim_\C \left( M_\rel (G)_{(p,w)} \otimes \dfrac{\O_{r+k,w}}{\m_{r+k,w}} \right)= \\&= \dim_\C \dfrac{\O_{B_\epsilon, p}}{J(g_w)} = \mu(g_w, p). 
    \end{align*}
    Hence, by Siersma's Theorem, it follows that $\sum_{p\in B_\epsilon / Y_w} \mu(g_w,p) = b_n(Y_w \cap B_\epsilon)$ is the $n$-th Betti number of $Y_w\cap B_\epsilon$, which is the image Milnor number $\mu_I(X,f)$. The result follows.
\end{proof}
With this result being proven, the following theorem is now an immediate consequence:
\begin{theorem}\label{theorem:MC} In the above conditions, the following statements are equivalent and imply the Mond conjecture for $(X,f)$:
\begin{enumerate}
    \item $\dim_\C M(g)=\mu_I(X,f)$,
    \item $M_\rel(G)$ is a Cohen-Macaulay $\O_{r+k}$-module of dimension $r+k$.  
\end{enumerate}
Furthermore, if $g$ is weighted homogeneous and satisfies the Mond conjecture, then the above assertions hold. 
\end{theorem}
    \begin{proof} This easily follows from \ref{app:samuel2}, which states that the Samuel multiplicity of a module $M$ is generally smaller than the dimension of $M/(\m \cdot M)$, with equality if and only if $M$ is Cohen-Macaulay with the same dimension as its base ring. Hence, $M_\rel (G)$ is Cohen-Macaulay of dimension $r+k$ if and only if
    \begin{align*}
        \mu_I(X,f)&=e\Big(\m_{r+k}, M_\rel (G)\Big) = \dim_\C \dfrac{M_\rel (G)}{\m_{r+k} M_\rel (G)}  = \\ 
        &= \dim_\C M_\rel (G)\otimes \dfrac{\O_{r+k}}{\m_{r+k}}=\dim_\C M(g).
    \end{align*}

    Therefore, if one of the items hold, it then follows that 
    \begin{align*}
        \mu_I(X,f)&=\dim_{\C} M(g) = \dim_\C K(g) + \codimAe (X,f) \geq \\ &\geq \codimAe (X,f),
    \end{align*}
    and hence, the generalised Mond conjecture holds for $(X,f)$. Moreover, if $g$ is weighted homogeneous, then $K(g)=0$. Thus, if the generalised Mond conjecture holds for $(X,f)$, then $\mu_I(X,f)=\dim_\C M(g)- \dim_\C K(g)=\dim_\C M(g)$. Hence, the above assertions hold. 
    \end{proof}
\begin{remark} Recent work by Nuño-Ballesteros and Giménez Conejero \cite{roberto2} has shown that the requirement for $M_\rel(G)$ to be $r+k$-dimensional can be eliminated from the second condition. Thus, we have that $\dim_\C M(g)=\mu_I(X,f)$ if and only if $M_\rel (G)$ is Cohen-Macaulay. Indeed, if the dimension of $M(g)$ is strictly less than $r+k$, then $\mu_I(X,f)$ is zero, which, according to \cite{juanjo-roberto}, implies that $f$ is a stable map-germ.
\end{remark}
This result shows that verifying the Cohen-Macaulay property of the module $M_\rel(G)$ is sufficient to establish the validity of the Mond conjecture for $f$. While it has been observed that this module is Cohen-Macaulay in every computed example, providing a rigorous proof of this statement remains an open question.

\subsection{Proof of the Mond conjecture for \textit{n\,=\,2} and mappings defined on \icis}\label{section:5_4}
In this section, we will achieve our second and final objective in the project, which is to show the validity of the generalised Mond conjecture for $n=2$ by employing the results from the previous section.\\

Appendix \ref{appendix:3} includes the properties of modules and the necessary notions that are employed throughout this section. Therefore, we suggest that readers first review this appendix to gain a comprehensive understanding of the section.\\

Building on the main result of the previous section, to establish the Mond conjecture for mappings $f:(X,0)\rightarrow (\C^3,0)$ where $(X,0)$ is a 2-dimensional \textsc{icis,} it suffices to check that $M_\rel (G)$ is a Cohen-Macaulay module of dimension $r+k$. Notice first that the dimension of $M_\rel (G)$ is $\leq r+k$, due to the fact that 
\begin{equation*}
    \dim_\C M_\rel (G)\otimes \dfrac{\O_{n+1+r+k}}{\m_{r+k}\cdot \O_{n+1+r+k}} =\dim_\C M(g) <+\infty.
\end{equation*}
Therefore, it is enough to show that $\depth M_\rel (G) \geq r+k$. To achieve this, we apply the depth lemma \ref{app:depthlemma} to the exact sequence
\begin{equation*}
    0\rightarrow M_\rel(G) \rightarrow \dfrac{\F_1(F)}{J_y(G)}\rightarrow \dfrac{C(F)}{J_{y,z}(G)\cdot \O_{n+k+r}} \rightarrow 0
\end{equation*}
to obtain that
\begin{equation*}
    \depth M_\rel (G)\geq \min \left\{ \depth \dfrac{\F_1(F)}{J_y (G)}, \depth \dfrac{C(F)}{J_{y,z} (G)\cdot \O_{n+k+r}} + 1 \right\}.
\end{equation*}
Hence, both terms of the minimum have to be shown to be greater than or equal to $r+k$. For the latter, it has been checked in lemma \ref{lemma:C/J_CM} that this module is Cohen-Macaulay of dimension $n+r+k-2$ for all values of $n$. In particular, it follows that 
\begin{equation*}
    \depth \dfrac{C(F)}{J_y (G)\cdot \O_{n+r+k}} + 1 = n+r+k-1 \geq r+k
\end{equation*}
for every $n\geq 1$. Therefore, it is enough to verify that 
\begin{equation*}
    \depth \dfrac{\F_1(F)}{J_y (G)} \geq r+k.
\end{equation*}
Notice that, up to this point, the assumption $n=2$ has not been used yet. In general, the module $\F_1(F)/J_y (G)$ has dimension 
\begin{equation*}
    \dim \dfrac{\F_1(F)}{J_y (G)} = \max \left\{\dim M_\rel (G), \dim \dfrac{C(F)}{J_y (G)\cdot \O_{n+r+k,0}} \right\} =n+r+k-2,
\end{equation*}
due to the fact that $\dim M_\rel (G) \leq r+k \leq n+r+k-2$ provided $n\geq 2$.
Therefore, it is not expected that $M_\rel(G)$ will be Cohen-Macaulay for every $n\geq 2$. The only case in which we can expect this is when $n=2$, since, in this case, its dimension is precisely $r+k$. To verify this claim, we make use of Pellikaan's theorem:
\begin{theorem}[Pellikaan, Section 3 of \cite{Pellikaan}] If $J\subset F\subset R$ are ideals of $R$ where $J$ is generated by $m$ elements, $\grade (F/J)\geq m$ and $\pd\, (R/F)=2$, then $F/J$ is a perfect module and grade$\,(F/J)=m$.
\end{theorem}
This result plays a crucial role in showing that the module is indeed Cohen-Macaulay, as it is proven in the following proposition:
\begin{theorem}\label{theorem:n=2} If $n=2$, then $\F_1(F)/J_y (G)$ is a Cohen-Macaulay module. 
\end{theorem}
\begin{proof} We follow the notation of the previous result, taking $R=\O_{n+1+r+k}$, $F=\F_1(F)$ and $J=J_y (G)$. Notice that the quotient ring $R/F=\O_{n+1+r+k}/\F_1(F)$ is Cohen-Macaulay of dimension $n+r+k-1$ (see remark \ref{remark:F1_determinantal}). Hence, by the Auslander–Buchsbaum formula,
\begin{align*}
    \pd\, (R/F) &= \depth R - \depth R/F = \dim R - \dim R/F = \\&= (n+1+r+k) - (n+r+k-1) = 2. 
\end{align*}
Moreover, $J=J_y (G)$ is clearly generated by $n+1$ elements, namely the partial derivatives of $G$ with respect to the variables $y_1, \ldots, y_{n+1}$. Furthermore, 
\begin{align*}
    \grade (F/J) &= \depth (\Ann (F/J), \O_{n+1+r})= \height \Ann (F/J) = \\
    &= \dim \O_{n+1+r}-\dim F/J= n+1+r-(n+r-2)=3.
\end{align*}
Since $n+1=3$ is our case, Pellikaan's theorem states that $F/J$ is a perfect $\O_{n+1+r+k}$-module. Furthermore, since $\O_{n+1+r+k}$ is a local Cohen-Macaulay ring, then $F/J$ is a Cohen-Macaulay module. 
\end{proof}
With this proof established, the main result of this section follows readily as an application of the results presented in the previous section.
\begin{theorem} Let $(X,0)$ be an \icis of dimension 2 and $f:(X,0)\rightarrow (\C^3,0)$ be an $\A$-finite mapping with image $(Y,0)$. Then, $\mu_I(f)\geq \codimAe (X,f)$, with equality if $(Y,0)$ is weighted homogeneous. In other words, the Mond conjecture holds for $f$.  
\end{theorem}
\begin{proof}
    The previous results imply that $M_\rel (G)$ is a Cohen-Macaulay module. Hence, Theorem \ref{theorem:MC} implies that the Mond conjecture holds for $f$. 
\end{proof}
In this setting, the image Milnor number therefore satisfies that $\mu_I(X,f)=\dim_\C M(g)$. This gives an operative method to compute this number, as the following example shows:
\begin{example} Let $(X,0)\subset (\C^3,0)$ be the hypersurface defined by $x^3+y^3-z^2=0$ and let $f:(X,0)\rightarrow (\C^3,0)$ be the $\A$-finite mapping $f(x,y,z)=(x,y,z^3+xz+y^2)$. In this case, $\hat{f}(x,y,z)=(x,y,z^3+xz+y^2, x^3+y^3-z^2)$ turns out to be a stable mapping. With Singular, we easily obtain that $\codimAe(X,f)=6$ and $\mu_I(X,f)=\dim_\C M(g)=6$. 
\end{example}

\newpage
\thispagestyle{empty}
\textcolor{white}{blank page}
\newpage
\appendix
\section{Prerequisites on commutative algebra for chapters 5 and 6}
In this chapter, we present a motley collection of definitions and results from the world of commutative algebra, that are mandatory to completely understand the developments of sections \ref{section:5_2}, \ref{section:5_3} and \ref{section:5_4}.
\subsection{Tensor product of modules}\label{appendix:1}
In this first section, we study the tensor product of modules and some relevant aspects concerning it. This first section is crucial for section \ref{section:5_2} regarding the specialisation of the module $M_\rel(G)$ to the original $M(g)$. We follow the clear exposition given in \cite{Eilenberg}, as well as some results that appear in \cite{Matsumura}. \\

Let $R$ be a ring (commutative with identity) and let $M,N$ be $R$-modules. A basic operation is the formation of the tensor product $M\otimes_R N$. In order to properly define it, let us consider the free group $F$ generated by pairs $(m,n)$ with $m\in M$ and $n\in N$, and let $R$ be the subgroup generated by the elements of the form 
\begin{align*}
    &(m+m',n)-(m,n)-(m',n), \\
    &(m,n+n')-(m,n)-(m,n'), \\
    &(m,\lambda n )- (\lambda m, n),
\end{align*}
for $m,m'\in M, n,n'\in N$ and $\lambda \in R$. We define the \textit{tensor product} $M\otimes_R N$ as the quotient $F/R$ regarded as an abelian group. The class of $(m,n)\in F$ in $M\otimes_R N$ is denoted as $m\otimes_R n$, or just $m\otimes n$ when the underlying ring is clear from the context. We then have the formal rules 
\begin{align*}
    &(m+m')\otimes n = m\otimes n + m'\otimes n, \\
    &m\otimes (n+n')=m\otimes n + m\otimes n', \\
    &(\lambda m) \otimes n = m \otimes (\lambda n),
\end{align*}
for $m,m'\in M, n,n'\in N$ and $\lambda \in R$. The natural action $\lambda (m\otimes n) := (\lambda m)\otimes n$ is therefore a well-defined multiplication $(R, M\otimes_R N)\rightarrow M\otimes_R N$ that induces an $R$-module structure in the tensor product $M\otimes_R N$. 

\begin{remark} Notice that the correspondence $\phi:M\times N \rightarrow M\otimes_R N$ given by $\phi(m,n)=m\otimes_R n$ is bilinear and satisfies $\phi(\lambda m, n )=\phi(m, \lambda n)$ for each $m\in M, n\in N$ and $\lambda \in R$. In fact, it can be shown that any other $R$-module $D$ with a bilinear mapping $f:M\times N\rightarrow D$ satisfying that $f(\lambda m, n)=f(m, \lambda n)$ admits a unique factorisation $f=g\circ \phi$ where $g:M\otimes_R N\rightarrow D$ is a homomorphism. This is usually shortened as
\begin{equation*}
    \begin{tikzcd}
	{M\times N} & {M\otimes_R N} \\
	& D
	\arrow["\exists !\, g", dashed, from=1-2, to=2-2]
	\arrow["f"', from=1-1, to=2-2]
	\arrow["\phi", from=1-1, to=1-2]
\end{tikzcd}
\end{equation*}
In fact, this gives an equivalent definition of the tensor product through this universal property. 
\end{remark}
A straightforward consequence of homomorphisms between tensor products is the following:
\begin{proposition}
    For every $R$-module homomorphisms $\phi:M\rightarrow M'$ and $\psi:N\rightarrow N'$, the correspondence $\phi\otimes \psi :M\otimes N\rightarrow M'\otimes N'$ given by $(\phi\otimes \psi)(m\otimes n)=\phi(m)\otimes \psi(n)$ for $m\in M$ and $n\in N$ is a well-defined $R$-module homomorphism. 
\end{proposition}
This result can be stated in the language of cathegory theory just by saying that $\otimes_R$ is a functor.  \\

Another immediate properties of the tensor product operation is that the mapping $\lambda \otimes m\mapsto \lambda m$ gives an isomorphism $R\otimes_R M\cong M$ and that the mapping $m\otimes n \mapsto n\otimes m$ gives an isomorphism $M\otimes_R N\cong N\otimes_R M$. Furthermore:
\begin{proposition} If $I\subset R$ is an ideal and $M$ is an $R$-module, then $M/(I\cdot M) \cong M \otimes_R (R/I)$.
\end{proposition}
\begin{proof} Notice that the correspondence $\phi:M\rightarrow M\otimes_R (R/I)$ given by $\phi(m)=m\otimes [1]$ is a surjective $R$-homomorphism with kernel $I\cdot M$. Then, the first isomorphism theorem provides the desired isomorphism. 
\end{proof}
The module that appears in the previous proposition is commonly referred to as the \textit{fibre} of $M$ with respect to $I$. This is a useful tool to specialise modules, since the fibre $M /(I\cdot M)$ is an $R/I$-module that eliminates the information contained in $M$ from $I$. \\

Recall that a sequence of $R$-modules and $R$-homomorphisms 
\begin{equation*}
    M_n \overset{\phi_n}{\rightarrow}M_{n-1} \overset{\phi_{n-1}}{\rightarrow}\ldots \overset{\phi_1}{\rightarrow}M_0
\end{equation*}
is said to be \textit{exact} whenever $\ker \phi_{i-1}= \im \phi_i $ for every $i\in \{1, \ldots, n\}$. We frequently consider \textit{short exact} sequences, which have the form
\begin{equation*}
    0 \rightarrow M'\overset{\phi}{\rightarrow} M \overset{\psi}{\rightarrow} M''\rightarrow 0
\end{equation*}
In this case, exactness is equivalent to have that $\phi$ is injective, $\psi$ is surjective and that $\im \phi = \ker \psi$. \\

If an exact sequence is tensored with a fixed module, then the sequence is no longer short exact. However, it turns out that 
\begin{proposition}
    If $0\rightarrow M'\rightarrow M\rightarrow M''\rightarrow 0$ is a short exact sequence of $R$-modules and $N$ is an $R$-module, then the induced sequence 
    \begin{equation}\label{eqn:rightexact}
        M'\otimes N \rightarrow M\otimes N\rightarrow M''\otimes N\rightarrow 0
    \end{equation}
    is exact. 
\end{proposition}
For a proof, see 4.5 of \cite{Eilenberg}. In the literature, this result is often shortened just by saying that the functor $\otimes$ is \textit{right exact}. \\

The kernel $K$ of the homomorphism $M'\otimes N \rightarrow M\otimes N$ from the exact sequence (\ref{eqn:rightexact}) is generally nonzero (in fact, if $K=0$, then the sequence is short exact). In the case that $M$ is a free module, it can be shown that $K$ depends only on $M''$ and $N$. We therefore define the \textit{torsion product} $\Tor^R_1(M'',N)$ to be this kernel. In this case, the sequence 
\begin{equation}\label{eqn:rightexact2}
        \Tor^R_1(M'',N)\rightarrow M'\otimes N \rightarrow M\otimes N\rightarrow M''\otimes N\rightarrow 0
    \end{equation}
    is exact, where the new homomorphism arising from the torsion group is just the inclusion map. Continuing this way, one can obtain an infinite exact sequence
    \begin{multline*}
        \ldots \rightarrow \Tor^R_{n+1}(M'',N) \rightarrow \Tor^R_{n}(M',N) \rightarrow \\ \rightarrow \Tor^R_{n}(M,N) \rightarrow \Tor^R_{n}(M'',N) \rightarrow \ldots
    \end{multline*}
    for $n\in \N$, which terminates in the original exact sequence given in (\ref{eqn:rightexact2}) provided we set $\Tor^R_0(M,N):=M\otimes_R N$. 
    \begin{definition}
        An $R$-module $M$ is $R$\textit{-flat} if, for every short exact sequence of $R$-modules $0\rightarrow N'\rightarrow N\rightarrow N''\rightarrow 0$, the associated sequence
        \begin{equation*}
            0 \rightarrow N'\otimes_R M \rightarrow N \otimes_R M \rightarrow N''\otimes_R M \rightarrow 0
        \end{equation*}
        is short exact. 
    \end{definition}
    The following proposition therefore becomes evident:
    \begin{proposition}
        An $R$-module $M$ is $R$-flat if and only if for every $R$-module $N$ one has that $\Tor^R_1(M,N)=0$. 
    \end{proposition}
    Flatness can be difficult to check in general. In our case, we are interested in the flatness of local rings. More precisely, let $R\subset S$ be local rings, in such a way that $S$ is an $R$-module, and let $\m$ be the maximal ideal of $R$. Denote by $F:=S/\m S \cong S\otimes_R (R/\m)$ the fibre of $S$ over $\m$. It can be checked that, if $S$ is $R$-flat, then 
    \begin{equation*}
        \dim S = \dim R + \dim F,
    \end{equation*}
    see Theorem 15.1 of \cite{Matsumura}. As the following result shows, a partial converse result holds:
    \begin{theorem}
        With the above notation, and if $R$ is a regular local ring, $S$ is a Cohen-Macaulay ring, and $\dim S = \dim R + \dim F$, then $S$ is $R$-flat. 
    \end{theorem}
    For a reference, see Theorem 23.1 of \cite{Matsumura}. \\

    An important result regarding the torsion groups is the following practical and elementary formula:
\begin{proposition}\label{prop:torsionIJ}
    If $I,J\subset R$ are ideals of the ring $R$, then $$\Tor^R_1\left(\dfrac{R}{I}, \dfrac{R}{J}\right)=\dfrac{I\cap J}{I\cdot J}.$$ 
\end{proposition}
For a proof, see \cite{Eilenberg} or \cite{Matsumura}.

\subsection{Samuel multiplicity of a module}\label{appendix:2}

In this section, we introduce a key concept that arises in one of the most relevant results of the project: the samuel multiplicity. In particular, we show in section \ref{section:5_3} that the image Milnor number is in fact the Samuel multiplicity of the module $M_\rel (G)$ over $\m_r$. In order to prove this result, it is crucial to know some technical aspects regarding this concept.  \\

Let $R$ be a Noetherian local ring of dimension $d$. We say that an ideal $\q$ of $R$ is $\m$\textit{-primary} provided $\m^k \subset \q \subset \m$ for some $k\geq 1$. As one might notice, this is equivalent to have that $\sqrt{\q}=\m$, which is in turn equivalent to have that $\length R/\q < \infty$ by proposition C.4 of \cite{juanjo}. In case that $\q$ is generated by $d$ elements, we say that $\q$ is a \textit{parameter ideal} of $R$. 
\begin{definition} Let $M$ be a finitely generated $R$-module, and let $\q$ be a parameter ideal of $R$. The \textit{Samuel function} of $M$ with respect to $\q$ is given by 
\begin{equation*}
    \chi_M^\q (t)= \length \left( \dfrac{M}{\q^{t+1}\cdot M} \right).
\end{equation*}
\end{definition}
It can be shown that this function is, in fact, a polynomial in $\Q[t]$ of degree $\dim M\leq d$ for large enough values of $t$. For a proof, we refer the reader to \cite{Matsumura}, pages 97 and 98. Furthermore, $\chi_M^\q$ can only take integer values. Hence, it can be shown through an induction on $d$ that 
\begin{equation*}
    \chi_M^\q (t) = \dfrac{e}{d!} t^d + \text{ terms of lower order,}
\end{equation*}
for some prescribed non-negative integer $e$.
\begin{definition}
    The number $e$ appearing in the previous equation is called the \textit{multiplicity} of $M$ with respect to $\q$, and it is denoted as $e(\q, M)$. In the case that $\q=\m$, we commonly denote $e(\q, M)=e(M)$ just as the multiplicity of $M$. 
\end{definition}
For the sake of brevity, we only present the two main results for which we are concerned. For a more detailed exposition, we refer the reader to section C6 of \cite{juanjo}. 
\begin{proposition}\label{app:samuel2}[C15 of \cite{juanjo}] For any finitely generated $\O_r$-module $M$, we have that 
\begin{equation*}
    e(M)\leq \dim_{\C} \dfrac{M}{\m_r \cdot M},
\end{equation*}
with equality if and only if $M$ is Cohen-Macaulay of dimension $r$. 
\end{proposition}
We say that a sheaf of modules $\M$ represents $M$, or that $\M$ is a representative of $M$, whenever $M$ is the stalk of $\M$ at the origin.
\begin{theorem}\label{app:samuel1}[E15 of \cite{juanjo}] Let $M$ be a finitely generated $\O_{X,x_0}$-module, let $f:(X,x_0)\rightarrow (\C^r, 0)$ be an analytic map such that $\dim_\C M/f^*\m_r M<\infty$ and let $\M$ be a representative of $M$. Then, there exists a representative $f:V\rightarrow W$ such that, for every $u\in W$, 
\begin{equation*}
    e(f^*\m_r, M)=\sum_{f(x)=u} e\Big(f^* \m_{r,u}, \M_{x}\Big).
\end{equation*} 
\end{theorem}
\subsection{Properties of modules}\label{appendix:3}
In this section, we study some properties of modules to measure its behaviour that are required in the development of sections \ref{section:5_2} and \ref{section:5_4}. First, we comment on some fact regarding the Cohen-Macaulay property and the depth of a module. Lastly, an analogous concept of Cohen-Macaulayness is defined, namely the concept of a perfect module. Its main properties are studied, as well as its relationship with the former concept. \\

The following result regarding the depth with exact sequences turns out to be crucial to determine the Cohen-Macaulay property of $M_\rel (G)$ for $n=2$ in chapter 6, namely the depth lemma:
\begin{lemma}[Depth lemma]\label{app:depthlemma} If $0\rightarrow M' \rightarrow M\rightarrow M''\rightarrow 0$ is a short exact sequence of nonzero finite $R$-modules, then
\begin{equation*}
    \depth M' \geq \min \{\depth M, \depth M'' + 1\}.
\end{equation*}
\end{lemma}
For a reference of these results, see the references \cite{Matsumura} or \cite{Eilenberg}. \\

In what follows, we define the notion of projective dimension. First, we recall that an $R$-module $P$ is said to be \textit{projective} whenever every homomorphism of $P$ into a quotient $N/N'$ admits a factorisation $P\rightarrow N \rightarrow N/N'$. Now, if $M$ is an $R$-module, then an exact sequence
\begin{equation}\label{eq:pd}
    \ldots \rightarrow P_n \rightarrow P_{n-1}\rightarrow P_1 \rightarrow P_0 \rightarrow M \rightarrow 0
\end{equation}
is called a \textit{projective resolution} of $M$ if each $P_i$ is projective, for $i\in \{0, 1, 2, \ldots \}$. If there exists a natural number $l$ such that, for each $i>l$ we have that $P_i=0$, and $P_{l}\neq 0$, we say that the length of the sequence is $l$. In this case, we call \textit{projective dimension} of $M$, and we denote it by $\pd_R(M)$, to the smallest $l$ among all possible projective resolutions (if $M$ has no finite projective resolutions, we then set $\pd_R(M)=\infty$). Whenever the base ring $R$ is clear from the context, we just write $\pd_R(M)=\pd\,(M)$. \\

A natural consequence is that the projective dimension of a projective module equals $0$. In fact, it is straightforward to notice that the projective modules are completely characterised as the modules with $\pd_R(M)=0$. \\

The projective dimension is tightly related with the notion of depth, as Auslander and Buchsbaum's formula states:
\begin{theorem}[Auslander-Buchsbaum, Theorem 19.1 of \cite{Matsumura}] Let $R$ be a Noetherian local ring and $M\neq 0$ be a finite $R$-module. Suppose that $\pd_R(M)<\infty$. Then,
\begin{equation*}
    \pd_R(M) + \depth M = \depth R. 
\end{equation*}
\end{theorem}

An interesting fact to take into account is that tensoring with a module $N$ in (\ref{eq:pd}) gives a sequence
\begin{equation*}
    \ldots P_n\otimes N \rightarrow \ldots \rightarrow P_0\otimes N,
\end{equation*}
which may not be exact, but which is a complex. It turns out that its $n$-th homology group is precisely $\Tor_n^R(M,N)$. \\

Another functor of at least as great importance as the tensor product is given by the group $\Hom_R(M,N)$ of $R$-homomorphisms $M\rightarrow N$. This functor can be shown to be contravariant in $M$, covariant in $N$, and left exact in the sense that, when applied to a short exact sequence $0\rightarrow M'\rightarrow M \rightarrow M'' \rightarrow 0$, we get that 
\begin{equation}\label{eq:hom}
    0  \rightarrow \Hom_R(M'',N) \rightarrow \Hom_R(M,N) \rightarrow \Hom_R(M',N) 
\end{equation}
is exact. In the same way that we defined $\Tor_n^R$, we define the \textit{ext} groups as the ones that fit into the exact sequence
\begin{multline*}
    \ldots \rightarrow \Ext_R^n(M'',N) \rightarrow \Ext_R^n(M,N) \rightarrow \\ \rightarrow \Ext_R^n(M',N) \rightarrow \Ext_R^{n+1}(M'',N)\rightarrow \ldots
\end{multline*}
which is a continuation of (\ref{eq:hom}) provided that we set $\Ext_R^0(M,N)=\Hom_R(M,N)$. These properties together with the property that
\begin{equation*}
    \Ext_R^n (P,N)=0 \text{ if }n>0 \text{ and }P \text{ is projective}
\end{equation*}
and the usual formal properties of functors suffice as an axiomatic description of the functors $\Ext_R^n(M,N)$. \\

Having this concept defined, one can now introuce the notion of \textit{grade} of a module, which was introduced first by Rees. If $M\neq 0$ is a finite $R$-module, we set
\begin{equation*}
    \grade (M) = \inf \{ n : \Ext_R^n (M,R)\neq 0\}.
\end{equation*}
Moreover, for a proper ideal $I$ of $R$ we call $\grade (R/I)$ the grade of the ideal $I$, and we just write $\grade (J)$. In virtue of Theorem 16.7 of \cite{Matsumura}, it turns out that $\grade (M) = \depth (\Ann (M), R)$, where $\Ann (M)= \{a\in R: a\cdot m =0 \, \forall m\in M\}$. Moreover, if $g=\grade (M)$, then $\Ext^g_R(M,R)\neq 0$, so that 
\begin{equation*}
    \grade (M) \leq \pd (M).
\end{equation*}
In a parallel way we did in the definition of Cohen-Macaulay modules, we say that $M$ is a \textit{perfect} module whenever equality holds in the previous equation, that is, if $\grade (M)=\pd (M)$. \\

As the definition may suggest, perfect modules and Cohen-Macaulay modules are extremely similar definitions. Indeed:

\begin{proposition} If $R$ is a local Cohen-Macaulay ring and $M$ is a perfect module, then $M$ is Cohen-Macaulay. 
\end{proposition}
\begin{proof} Notice that, if $I=\Ann (M)$, then
\begin{align*}
    \pd (M) &= \grade (M) = \depth (I,R) =  \height (I) = \\
    &= \dim R - \dim (R/I) =\dim R- \dim M,
\end{align*}
where Theorem 17.4 of \cite{Matsumura} is applied. Taking this into account, it follows by the Auslander-Buchbaum formula that $\pd (M)= \depth R-\depth M=\dim R - \depth M$, due to the fact that $R$ is Cohen-Macaulay. Therefore, $\dim M= \depth M$, and the claim follows.
\end{proof}
In fact, it can be shown that, if $R$ is a regular ring, then $M$ is a perfect $R$-module if and only if it is Cohen-Macaulay. Therefore, the differences between these modules can only be noticed over nonregular rings. \\

The main result that allows us to prove the conjecture for $n=2$ is the following result from Pellikaan, which states that:
\begin{theorem}[Pellikaan, Section 3 of \cite{Pellikaan}] If $J\subset F\subset R$ are ideals of $R$ where $J$ is generated by $m$ elements, $\grade (F/J)\geq m$ and $\pd\, (R/F)=2$, then $F/J$ is a perfect module and grade$\,(F/J)=m$.
\end{theorem}

In this case, our intention is to verify that the module $F/J$ is Cohen-Macaulay, and this is achieved with a combination of Pellikaan's theorem with the previous definition. For more details, see the Theorem \ref{theorem:n=2}. 

\newpage
\addcontentsline{toc}{section}{References}{\protect\numberline{}}
\nocite{*}
\printbibliography

@article{vanishingcycles,
  author = {D. Mond},
  title = {Vanishing cycles for analytic maps},
  journal = {Singularity Theory and Its Applications, Part I},
  year = 1991,
  volume = {1462},
  pages = {221–234},
}

@article{BobadillaNunoPenafort,
  author = {J. Fernández de Bobadilla and J.J. Nuño-Ballesteros and G. Peñafort-Sanchis},
  title = {A Jacobian module for disentanglements and applications to Mond’s conjecture},
  journal = {Rev Mat Complut},
  year = 2019,
  volume = {32},
  pages = {395–418},
}

@article{Pellikaan,
  author = {R. Pellikaan},
  title = {Projective resolutions of the quotient of two ideals},
  journal = {Indagationes Mathematicae (Proceedings)},
  year = 1988,
  volume = {91},
  pages = {65-84},
}

@article{roberto,
  author = {R. Giménez Conejero and J.J. Nuño-Ballesteros},
  title = {Singularities of mappings on ICIS and applications to Whitney equisingularity},
  journal = {Advances in Mathematics},
  year = 2022,
  volume = {408},
}

@BOOK{juanjo,
AUTHOR={D. Mond and J.J. Nuño-Ballesteros},
TITLE={Singularities of mappings - The local behaviour of smooth and complex analytic mappings},
PUBLISHER={Grundlehren der Mathematischen Wissenschaften, 357. Springer, Cham},
YEAR={2020}
}

@article{Piene,
  title={Ideals associated to a desingularitation},
  author={R. Piene},
  journal={Proceedings of the SummerMeeting Copenhagen
1978. Lecture Notes in Mathematics},
  volume={732},
  number={},
  pages={503–517},
  year={1979},
  publisher={Springer, Berlin}
}

@article{Mond-Pellikaan,
  title={Fitting ideals and multiple points of analytic mappings},
  author={D. Mond and R. Pellikaan},
  journal={Algebraic Geometry
and Complex Analysis (Pátzcuaro, 1987). Lecture Notes in Mathematics},
  volume={1414},
  number={},
  pages={107–161},
  year={1989},
  publisher={Springer, Berlin}
}

@article{Mond-Montaldi,
  title={Deformations of maps on complete intersections, Damon’s KV-equivalence and bifurcations},
  author={D. Mond and J. Montaldi},
  journal={Singularities, Lille, 1991, in: London Math. Soc. Lecture Note Ser.},
  volume={201},
  number={},
  pages={263–284},
  year={1994},
  publisher={ Cambridge Univ. Press, Cambridge}
}

@article{juanjo-roberto,
    author = {R. Giménez Conejero and   J.J. Nuño-Ballesteros},
    title = {The Image Milnor Number And Excellent Unfoldings},
    journal = {The Quarterly Journal of Mathematics},
    volume = {73},
    number = {1},
    pages = {45-63},
    year = {2021},
    month = {03},
    issn = {0033-5606},
}

@article{juanjo-henrique,
author = {D.A. Henrique and J.J. Nuño-Ballesteros},
title = {Mond's conjecture for maps between curves},
journal = {Mathematische Nachrichten},
volume = {290},
number = {17-18},
pages = {2845-2857},
year = {2017}
}

@BOOK{Eilenberg,
AUTHOR={H. Cartan and S. Eilenberg},
TITLE={Homological algebra},
PUBLISHER={Princeton University Press},
YEAR={1956}
}

@BOOK{Matsumura,
AUTHOR={H. Matsumura},
TITLE={Commutative algebra},
PUBLISHER={Benjamin},
YEAR={1970}
}

@article{Siersma,
  title={Vanishing cycles and special fibres},
  author={D. Siersma},
  journal={Singularity Theory and Its Applications, Part I
(Coventry, 1988/1989). Lecture Notes in Mathematics,},
  volume={1462},
  number={},
  pages={292–301},
  year={1991},
  publisher={Springer, Berlin}
}

@BOOK{Looijenga,
AUTHOR={E. J. N. Looijenga},
TITLE={Isolated Singular Points on Complete Intersections},
PUBLISHER={Cambridge University Press},
YEAR={1984}
}

@BOOK{icis,
AUTHOR={T. Gaffney and J. J. Nuño-Ballesteros and B. Oréfice-Okamoto and M. A. S. Ruas & R. Wik-Atique},
TITLE={Complete intersections with isolated singularities (ICIS). Algebraic methods and singularities},
PUBLISHER={CIMPA research school, Singularities and its applications. ICMC–USP São Carlos, Brazil},
YEAR={2022}
}

@article{roberto2,
  author = {R. Giménez Conejero and J.J. Nuño-Ballesteros},
  title = {A weak version of Mond's conjecture},
  journal = {arXiv:2207.01735, preprint},
  year = 2022,
}
\end{document}